\documentclass[dvips]{article}
\usepackage{amsmath,amsfonts,graphicx,amssymb,framed,relsize}
\usepackage{setspace}
\usepackage{graphicx}
\usepackage[english]{babel}
\usepackage{chngcntr,caption}
\usepackage{float}
\usepackage[affil-it]{authblk}
\restylefloat{table}
\begin{document}
\singlespacing
\title{The Noise Handling Properties of the Talbot Algorithm for Numerically Inverting the Laplace Transform}
\author{Colin L. Defreitas \&  Steve J. Kane}
\affil{School of Physics Astronomy and  Mathematics.}
\affil{University of Hertfordshire.}
\maketitle
\section{Abstract}
{\footnotesize This paper examines the noise handling properties of three of the most widely used algorithms for numerically inverting the Laplace Transform. After examining the genesis of the algorithms,
the regularization properties are evaluated through a series of standard test functions in which noise is added to the inverse transform.
Comparisons are then made with the exact data. Our main finding is that the Talbot inversion algorithm is very good at handling noisy data and performs much better than the Fourier Series and Stehfest numerical inversion schemes as outlined in this paper. This offers a considerable advantage for it's use in inverting the Laplace Transform when seeking numerical solutions to time dependent differential equations.}\\

\section{The Laplace Transform}
The Laplace Transform is an integral transform defined as follows:\\
Let $f(t)$ be defined for $t>0$, then the Laplace transform of $f(t)$ is given by:
\begin {equation}\label{eq:1}
\mathcal{L}\{f(t)\}=\int_0^\infty f(t)e^{-st}\;\ dt 
\end{equation} \label{eq:1a}
with $\mathcal{L} \{f(t)\}$ denoted as $ F(s)$. The Laplace transform can be shown to exist for any function which can be integrated over any finite interval $0<t<l$
for $l>0$, and for which $f(t)$ is of exponential order, i.e.
\begin{equation}
\mid f(t)\mid <Me^{at}
\end{equation}
as $t\rightarrow \infty$,
where $M$ and $a$ are some small real positive numbers.\\
Analytically the inverse Laplace transform is usually obtained using the techniques of contour integration with the resulting set of  standard transforms presented in tables e.g. Spiegel (1965). However, using the Laplace Transform to obtain solutions of differential equations can lead to solutions in the Laplace domain which are not easily invertible
to the real domain by analytical means. Thus Numerical inversions are used to convert the solution from the Laplace domain to the real domain.\\
For example one of the principal reasons for employing the Laplace Transform in finding solutions to diffusion type problems is that in general, their numerical solutions invariably employ finite difference (FD) schemes involving grid representations which for convergence to a desired solution can preclude the use of large time-steps as these can seriously affect the accuracy, and for explicit schemes, the stability of the solution.\\
For these reasons FD schemes may require many very small time steps to reach the desired solution at any given time. This can lead to increases in the round-off error 
inevitably generated in these calculations.\\
Moreover as Mordis (1994) points out, ``since these calculations essentially involve solving a large system of linear equations, FD methods may require several hundred matrix inversions. But the Laplace Transform Finite Difference Method, LTFDM requires only one time-step and no more than 6 to 10 matrix inversions".

The problems associated with the Finite Difference method become more acute in solving non-linear diffusion problems where now further iterations are required at each time-step. (Mordis(1984)).
Thus by employing the LTFDM  we can calculate the dependent variable for any time step desired without loss of accuracy or stability and without a large number of intermediate time steps.

\section{The Inverse Laplace Transform Perturbation and Precision }

The recovery of the function $f(t)$ is via the inverse Laplace transform which is most commonly defined by the Bromwich integral formula 
\begin{equation}\label{eq:2}
\mathcal{L}^{-1}\{F(s)\}=f(t)=\frac{1}{2\pi i}\int^{u +i\infty }_{u -i\infty }F(s)\; e^{st}\; ds,
\end{equation}
for some $u \in  \mathbb{R}$. (Spigel 1965)

The the choice of $s$ in (1) and so in (3) is not an arbitrary one. To demonstrate this we elaborate on some of the points made by Bellman and Roth (1984). We consider again the Laplace Transform
in (\ref{eq:1})
where we allow $s$ to be any complex number. As $s  \in \mathbb{C}$ then $F(s)$ must also be a function of a complex variable. From a numerical standpoint, the limits of integration of the Laplace Transform can over complicate the evaluation of the integral. This is  because the integral is evaluated over the semi infinite domain. To overcome this Bellman applies the transformation $r=e^{-t}$  and the integral becomes,
$\mathbb{c}$
\begin{equation}\label{eq:4} 
-\int^0_1 f(-\ln r) \;\frac{r^s}{r} dr =\int^1_0 \;\ {r^{s-1}}\;\ g(r)\;\  dr 
\end{equation}
where $g(r)=f(-\ln r)$.\\
As $s  \in \mathbb{C}$ then with $s = u +iv$ we have
\begin {equation}\label{eq:5}
r^{s-1}=r^{u-1+iv}=r^{u-1}e^{iv\ln r}.
\end{equation}
Now $F(s)=a(u,v)+ib(u,v)$ thus
\begin{equation}\label{eq:6}
a(u,v)=\int^{1}_{0}r^{u-1}\cos (v\ln r) g(r)dr
\end{equation}
\begin{equation}\label{eq:7}
b(u,v)=\int^{1}_{0}r^{u-1}\sin (v\ln r)g(r)dr.
\end{equation}
We now have two real integrals to evaluate. More importantly ``if $u<1$ there will be a singularity at $r=0$ which makes numerical computation difficult if not impossible.
In addition if $v\neq 0$ aside from the rapid oscillations at the origin the integral is unbounded for $u<1$"' . (Bellman and Roth, 1984).\\
To eliminate these problems  Bellman and Roth  choose $s$ to lie on the positive real axis greater than unity. In choosing $s$ to lie on the positive real axis we are treating the solution of (3) as a positive real integral equation. The problem here is that the inverse problem is known to be ill-posed meaning that small changes in the values of $F(s)$ can lead to large errors in the values for $f(t)$. (Davies and Martin, 2001).

This suggests that  LTFD methods which rely on such inversion methods can be highly sensitive to the inevitable noisy data that arises in their computation via truncation and round off error, a process which is exacerbated in non-linear schemes. While as Abate and Valko (2004) have shown  this to some extent can be curtailed by working in a multi-precision environment, as we show below a small amount of noise in the data can still enormously perturb the solution. When this is the case it becomes difficult for unlimited precision to aid in the convergence of the algorithm to the correct solution.\\

\section{The Algorithms}
There are over 100 algorithms available for inverting the Laplace Transform with numerous comparative studies, examples include: Cost (1964), Davies and Martin (1979), Duffy (1993), Narayanan and Beskos (1982) and Cohen (2007). However for the purposes of this investigation we follow Abate and Valko (2004) and apply our tests using ``Those algorithms that have passed the test of time". These fall into four groups,
\begin{flushleft}
	(1)   Fourier series expansion.\\
	(2)  Combination of Gaver Functionals.\\
	(3) Laguerre function expansion.\\
	(4)  Deformation of the Bromwich contour.\\
\end{flushleft}

Derivations of particular versions of these algorithms are given in the next section. However for now we do not run our tests using the Laguerre function expansion. While we do intend to investigate this method later on in our work our choices in this work have been made based on the ease of implementation of the inversion method an issue connected to parameter choice and control. The Laguerre method requires more than two parameters to effectively compute the desired transform while the other three methods can perform reasonably well when defined using just the one parameter.

\subsection{The Fourier Series Method}
In their survey of algorithms for inverting the Laplace Transform, Davies and Martin (1979)  note that the Fourier series method without accelerated convergence gives good accuracy on a wide variety of functions. Since the Laplace Transform is closely related to the Fourier transform it is not surprising that inversion methods based on a Fourier series expansion would yield accurate results. In fact the two sided Laplace transform can be derived from the Fourier transform in the following way.
We can define the Fourier transform as
\begin{equation}\label{eq:32}
\mathcal{F} \{f(t)\} = \int^\infty_{-\infty} f(t) \;\ e^{-2\pi i\nu t}\;\ dt
\end{equation}
then letting $v = 2\pi \nu $ we have 
\begin{equation}\label{eq:3a}
\mathcal{F} \{f(t)\} = \int^\infty_{-\infty} f(t) \;\ e^{-iv t}\;\ dt
\end{equation}
providing $f(t)$ is an absolutely integrable function, i.e. 
\begin{equation}\label{eq:33}
\int^\infty_{-\infty} |f(t)| \;\ dt \;\ < \infty.
\end{equation}
As many functions do not satisfy the condition in (\ref{eq:33}), $f(t)$ is multiplied by the exponential dampening factor $e^{-u t}$ thus
\begin{equation}\label{eq:34}
\mathcal{F} \{f(t)e^{-u t}\} = \int^\infty_{-\infty} f(t) \;\ e^{-iv t}e^{-u t}\;\ dt
\end{equation}
and letting $s=u + iv$  we obtain the two sided Laplace Transform of $f(t)$ as
\begin{equation}\label{eq:35}
\mathcal{F} \{f(t)e^{-u t}\} = \mathcal{L} \{f(t)\}=\int^ \infty_{-\infty} e^{-st} f(t) \;\ dt.
\end{equation}
LePage (1961) notes that the integral in (\ref{eq:35}) can be written in two parts as follows:
\begin{equation}\label{eq:35a}
\int^ \infty_{-\infty} e^{-st} f(t) \;\ dt =\int^ 0_{-\infty} e^{-st} f(t) \;\ dt +\int^ {\infty}_0 e^{-st} f(t) \;\ dt
\end{equation}
\begin{equation}\label{eq:35a}
=\int^ {\infty}_0 e^{-st} f(-t) \;\ dt  +\int^ {\infty}_0 e^{-st} f(t) \;\ dt
\end{equation}
The second term in the above expression is refereed to as the one sided Laplace transform or simply the Laplace transform.
Thus $s$ is defined as a complex variable in the definition of the Laplace Transform.\\

As before the inverse Laplace transform is given as:
\begin{equation}\label{eq:38}
f(t) = \frac{1}{2\pi i}\int^{u+ i\infty}_{u - i\infty}e^{st} F(s)\;\ ds.
\end{equation}
With $s = u +iv$ in (\ref{eq:38}) this leads to the result 
\begin{equation}\label{eq:39}
f(t) = \frac{2e^{ut}}{\pi}\int^\infty_0 [\mbox{Re} \{F(u+iv)\}\cos\ (v  t) - \mbox{Im }\{F(u+iv)\}\sin\ (v  t)]\;\ dv. 
\end{equation}
As Crump (1976) points out equations  (\ref{eq:1}) and (\ref{eq:39}) can be replaced by the cosine transform pair
\begin{equation}\label{eq:40}
\mbox{Re} \{F(u+iv)\} = \int^\infty_0 e^{-ut}\;\ f(t) \cos\ (v t)\;\ dt
\end{equation}
\begin{equation}\label{eq:41}
f(t) = \frac{2e^{ut}}{\pi}\int^\infty_0 \mbox{Re} \{F(u+iv)\} \cos (v t)\;\ dv
\end{equation}
or by the sine transform pair
\begin{equation}\label{eq:42}
\mbox{Re} \{F(u+iv)\} = -\int^\infty_0 e^{-ut}\;\ f(t) \sin\ (v t)\;\ dt
\end{equation}
\begin{equation}\label{eq:43}
f(t) = -\frac{2e^{ut}}{\pi}\int^\infty_0 \mbox{Im}  \{F(u+iv)\} \sin (v t)\;\ dv.
\end{equation}

Dunbar and Abate (1968) applied a trapezoid rule to (\ref{eq:41}) resulting in the Fourier series approximation,
\begin{equation}\label{eq:44}
f(t)\approx \frac{2e^{ut}}{T}\bigg[ \frac{1}{2}\;\ F(u) + \sum^\infty_{k=1}\mbox{Re} \bigg\{F\Big(u + \frac{k\pi i}{T}\Big)\bigg\} \cos\bigg(\frac{k\pi t}{T}\bigg)\bigg],
\end{equation}
where $f(t)$ is expanded in the interval $0 \leq t < T$.
For faster computation Simon and Stroot (1972) suggests that we let $T = 2t$ so we have
\begin{equation}\label{eq:45}
f(t)\approx \frac{e^{ut}}{t}\bigg[ \frac{1}{2}\;\ F(u) + \sum^\infty_{k=1}\mbox{Re} \bigg\{F\Big(u + \frac{k\pi i}{t}\Big)\bigg\} (-1)^k\bigg].
\end{equation}
This series can be summed much faster than (\ref {eq:44}) as there are no cosines to compute. (Crump,1976).
This algorithm is relatively easy to implement with $u$ being the only real varying parameter. It also has the advantage that there are  no deformed contours to consider so little prior knowledge of singularities and therefore of the nature of the solution to a particular numerical inversion are necessary allowing for a wide variety of applications.\\
However as Crump (1976)  pointed out for the the expression in (\ref{eq:45}) the transform $F(s)$ must now be computed for a different set of $s-$values for each distinct $t$. Since this type of application occurs often in practice  in which the numerical computations of $F(s)$ is itself quite time consuming  this may not be an economical inversion algorithm to use. These drawbacks  to some extent can be overcome by using the fast Fourier transform techniques (Cooley and Tukey 1965; Cooley et al.,1970).\\
Crump (1976) also  extends this method to one of faster convergence by making use of the the already computed imaginary parts. There are several other acceleration schemes for example those outlined by Cohen (2007). However these acceleration methods in general require the introduction of new parameters thus complicating the ease of implementation. 

\subsection{The Stehfest Algorithm}

In their survey Davies and Martin (1979) cite the Stehfest (1970) algorithm as providing accurate results on a variety of test functions. Since that time, this algorithm has become widely used for inverting the Laplace Transform being favored due its reported accuracy and ease of  implementation.\\
Below we give a brief overview of the evolution of the algorithm from a probability distribution function to the Gaver functional whose asymptotic expansion leads to an acceleration scheme which yields the algorithm in its most widely used form.\\
Gaver (1965) was investigating a method for obtaining numerical information on the time dependent behavior of stochastic processes which often arise in queuing theory. The investigation involved examining the properties of the three parameter class of density functions namely
\begin{equation}\label{eq:46}
p_{n,m}(a;t) = \frac{(n+m)!}{n!(m-1)!}(1 - e^{-at})^na e^{-mat}.
\end{equation}
After the binomial expansion of the term $(1 - e^{-at})^n$, Gaver went on to find the expectancy $E[f(T_{n,m})]$ where $T_{n,m}$ is the random variable with density 
(\ref{eq:46}). From this Gaver was able to express the inverse Laplace transform in terms of the functional
\begin{equation}\label{eq:47}
{f}_{n,m}(t) = \frac{\ln 2}{t}\frac{(n+m)!}{n!(m-1)!} \;\ \sum^n_{j=0}{n\choose k}(-1)^k F\bigg{(}(k+m)\frac{\ln 2}{t}\bigg{)}.
\end{equation}
With certain conditions on $n$ and $m$,  Gaver makes $n=m$ and expresses  (\ref{eq:47}) as
\begin{equation}\label{eq:48}
f_n(t)=\frac{\ln 2}{t}\frac{(2n)!}{n!(n-1)!} \sum^n_{k=0} {n\choose k} (-1)^k F\bigg{(}(k+n)\frac{\ln 2}{t}\bigg{)}
\end{equation}
While the expression in (\ref{eq:48}) can be used to successfully invert the Laplace Transform for a large class of functions, Gaver (1965) has shown that (\ref {eq:48}), with $a=\frac{\ln 2}{t}$ 
the Gaver functional has the asymptotic expansion
\begin{equation}\label{eq:49}
f_n(t) \approx f\bigg(\frac{\ln 2}{a}\bigg) + \frac{\alpha _1}{n}+\frac{\alpha _2}{n^2}+\frac{\alpha _3}{n^3} +...,
\end{equation} converging to the limit
\[ f\bigg(\frac{\ln 2}{a} \bigg)\] as $n\rightarrow\infty$.
(For the conditions on $m$ and $n$ and justification for the substitution for $a$ referred to above see (Gaver 1965)).
This asymptotic expansion provides scope for applying various acceleration techniques enabling a more viable application of the basic algorithm. Much of the literature eludes to to the fact that a Salzer (1956) acceleration  scheme is used on the Gaver functional in (\ref{eq:48}) which results in the Stehfest algorithm. In fact Stehfest's approach was a little more subtle than a direct application of the Salzer acceleration. They both however amount to the same thing.\\

\subsubsection{Using Salzer acceleration}

The Salzer acceleration scheme makes use of the ``Toeplitz limit theorem". Wimp (1981).  ``This concerns the convergence of a transformation of a sequence $\zeta _s$ where the $(n+1)$th member of the transformed sequence is a weighted mean of the first $(n+1)$ terms" 
\begin{equation}\label{eq:50}
\overline S_n =\sum^n_{k=0}\mu _{nk}.S_k.
\end{equation}
Here $\overline S_n $ is the transformed sequence and $S_k$ the original sequence. The Salzer means are given by 
\begin{equation}\label{eq:51}
\mu _{nk} = (-1)^{n+k}\;\ \frac{(1+k)^n}{n!} \;\ {n\choose k} .\;\
\end{equation}
For our purposes, $S_k = f_n(t)$ in (\ref{eq:48}). For the sake of compatibility with (\ref{eq:50}) we make the change $k \rightarrow i$ and $n\rightarrow k$ in (\ref{eq:48}). With this change of 
variables we also write
\[ \frac{(2k)!}{k!(k-1)!} = \frac{k(2k)!}{(k!k!)}\]
This allows the sum to be taken from $k=0$ to $n$ without $(0-1)!$ in the denominator of (\ref{eq:48}).
So with Salzer acceleration we can express our algorithm as
\begin{equation}\label{eq:52}
f_n(t)= \frac{\ln 2}{t}\sum^n_{k=0}(-1)^{n+k}\frac{(k+1)^n}{k!(n-k)!}\;\ \frac{k(2k)!}{k!k!}\sum^k_{i=0}\frac{k!}{i!(k-i)!}(-1)^i\;\ F\bigg\{\frac{(k+i)\ln 2}{t}\bigg\}.
\end{equation}

\subsubsection{Stehfest's approach}

For the purposes of following Stehfest's derivation it will be conveinent to rewrite (\ref{eq:48})  as
\begin{equation}\label{eq:53}
f_n(t) = {F}_n = \frac{(2n)!a}{n!(n-1)!} \;\ \sum^n_{j=0}{n\choose k}(-1)^k F\bigg{(}(k+n)a\bigg{)}
\end{equation}
with $a=\frac{\ln 2}{t}$.
Stehfest (1970) begins by  supposing  we have  $N$ values for $F[(k+n)a]$ with $F(a)$, $F(2a)$, $F(3a)$, ....$F(Na)$ for $N$ even. Using (\ref{eq:53}) we can then determine $\frac{N}{2}$ values $F_1,F_2,...,F_{N/2}$. Now  each of these $N/2$ values satisfy the asymptotic series in (\ref{eq:49}) with the same coefficients. 
This allows us to eliminate the first $N/2-1$ error terms.\\ We demonstrate the process for $N=6$, so we compute $3$ values. These have the following asymptotic expansions;
\[
{F}_1 = f\bigg(\frac{\ln 2}{a}\bigg) + \frac{\alpha _1}{1}+\frac{\alpha _2}{1^2}+\frac{\alpha _3}{1^3}+....
\] 
and as the $\alpha _j$'s are the same for each of the expansions then by using a suitable linear combination we can eliminate the first ($\frac{N}{2} -1$) error terms. That is
\begin{equation}\label{eq:54}
 f\bigg(\frac{\ln 2}{a}\bigg) = \sum^{\frac{N}{2}}_{n=1}\;\ a_n {F}_{(\frac{n}{2}+i-1)} + O\bigg(\frac{1}{N^{\frac{N}{2}}}\bigg)
\end{equation}
which may be achieved by selecting the coefficients to satisfy 
\begin{equation}\label{eq:55}
\sum^{\frac{N}{2}}_{n=1}\;\ a_n\frac{1}{(\frac{N}{2}+1-n)^k} =\delta_{k,0} \;\ \;\ \;\  k= 1,...,N/2 - 1.
\end{equation}
Hence we have 
\[
f\bigg(\frac{\ln 2}{a}\bigg) = a_1\bigg(f\bigg(\frac{\ln 2}{a}\bigg)+  \frac{\alpha _1}{1}+\frac{\alpha _2}{1^2}\bigg)
\]
\[
+ a_2\bigg(f\bigg(\frac{\ln 2}{a}\bigg)+  \frac{\alpha _1}{2}+\frac{\alpha _2}{2^2}\bigg)
\]
\begin{equation}\label{eq:56}
+ a_3\bigg(f\bigg(\frac{\ln 2}{a}\bigg)+  \frac{\alpha _1}{3}+\frac{\alpha _2}{3^2}\bigg)
\end{equation} 
thus
\[
a_1 + a_2 + a_3 =1
\]
\[
\alpha _1(a_1 + a_2 + a_3) =0
\]
\begin{equation}\label{eq:57}
\alpha _2(a_1 + a_2 + a_3) =0
\end{equation}
with  $\alpha _j\not = 0$ we have 
\[ a_1=\frac{1}{2},\;\ a_2=-4,\;\ a_3=\frac{9}{2}.
\]
Which are the same coefficients produced by the Salzer acceleration scheme used in (\ref{eq:50}) . In fact for any $n$, Stehfest generates the required coefficients using what is in effect a modified Salzer acceleration scheme. Giving
\begin{equation}\label{eq:58}
a_n = \frac{(-1)^{n-1}}{(\frac{N}{2})!}\;\ {\frac{N}{2}\choose n} \;\ n\bigg(\frac{N}{2}+1 -n)^{\frac{N}{2}-1}\bigg).
\end{equation}
Finally Stehfest substitutes these results into (\ref{eq:54}) and gets the inversion formula
\begin{equation}\label{eq:59}
f(t) \approx \frac{\ln 2}{t} \sum_{j=1}^N A_jF\bigg(\frac{j\ln 2}{t}\bigg)
\end{equation}
where 

for $N$ even.
\begin{equation}\label{eq:60}
A_j = (-1)^{\frac{N}{2}+j} \\= \sum^{min(j,\frac{N}{2})}_{k=\lfloor\frac{j+1}{2}\rfloor}\frac{k^{\frac{N}{2}}(2k)!}{(\frac{N}{2}-k)!k!(k-1)!(j-k)!(2k-j)!}.\\
\end{equation}
\subsection{The Talbot Algorithm.}
As we established in equations  (\ref{eq:32}) to (\ref{eq:35}) the Laplace transform  can be seen as a Fourier transform of the function
\begin{equation}\label{eq:61}
e^{-ut}f(t), \;\ t>0.
\end{equation}  
i.e.
\begin{equation}
{\cal F}\{e^{-ut}f(t)\}=\mathcal{L}\{f(t)\}={F(s)}
\end{equation}
Hence the Fourier transform inversion formula can be applied to recover the function thus:
\begin{equation}\label{eq:62}
{\cal F}^{-1}\;\ \{F(s)\} \;\  = \;\ e^{-ut}f(t) = \;\ \frac{1}{2\pi}\int^\infty_{-\infty}F(s)\;\ e^{ivt} \;\ dv
\end{equation}

and this equals
\begin{equation}\label{eq:63}
\frac{1}{2\pi i}\int^\infty_{-\infty} F(s) e^{st} \;\ dv
\end{equation}
 as $s = u + iv$ \;\ we have that \;\ $ds = \;\ idv$\;\ and so
 
 \begin{equation}\label{eq:64}
 f(t) \;\ = \frac{1}{2\pi i}\int^{u+i\infty}_{u-i\infty}F(s)\;e^{st}\;\ ds.
\end{equation}
 This result provides a direct means of obtaining the inverse Laplace transform. 
In practice the integral in (\ref{eq:64}) is evaluated using a  contour 
\begin{equation}\label{eq:65}
\frac{1}{2\pi i}\int_B e^{st}\;\ F(s)\;\ ds
\end{equation}
with $B$ here denoting the Bromwich contour. (Spigel, 1965). The contour is chosen so that it encloses all the possible singularities of $F(s)$.
The idea of the contour is introduced so that the  residue theorem can be used  to evaluate the integral.

Of course the residue theorem is used in the case of  an analytically invertible function. When $f(t)$ is to be calculated using numerical quadrature it may be more appropriate to devise a new contour. To ensure the convergence of (\ref{eq:65}) we may wish to control the growth of the magnitude of the integrand $e^{st}$ by moving the contour to the left so giving the real part of  $s$ a large negative component. (See Abate \& Valko (2003), Murli \& Rizzardi (1990) )\\
 But the deformed contour must not be allowed to pass through any singularities of $F(s)$. This is to ensure that the transform is analytic in the region to the right of $\textit B$.\\
\subsubsection{Derivation of the Fixed Talbot Contour.}

Talbot (1978) and Abate and Valko (2003) are used as the primary basis for extending the explanation of the derivation of the Talbot algorithm for inverting the Laplace Transform.\\
Abate and Valko begin with the Bromwich inversion integral along the Bromwich contour $B$  with the transformation

\begin{equation}\label{eq:66}
f\hat(s) = \frac{1}{s^  {\alpha} }, \;\ \alpha > 0.
\end{equation}
so $f(t)$ can be expressed as
\begin{equation}\label{eq:67}
f(t) = \frac{1}{2\pi i} \int_{B} e^{st} s^{-at} \;\ ds
\end{equation}
and thus 
\begin{equation}\label{eq:68}
f(t) = \frac{1}{2\pi i} \int_{B} e^{t(s - alog_e s)} \;\ ds
\end{equation}
with $a =\frac{\alpha}{t}$ in (\ref{eq:67}) and (\ref{eq:68}). As Abate and Valko (2003) point out numerically evaluating the  integral in (\ref{eq:68}) is difficult due to the oscillatory nature of the integrand.\\
However  this evaluation can be achieved by deforming the contour $B$ into a \textit{path of constant phase} thus eliminating the oscillations in the imaginary component. These paths of constant phase are also paths of steepest decent for the real part of the integrand. (See  Valko, 2003, Meissen, 2013, Bender and Orszag ,1978).\\
There are in general a number of contours for which the imaginary component remains constant so we choose one on which the real part attains a maximum on the interior (a saddle point) and this occurs at $g^{'}(s) = 0$ at some point on the contour. At these saddle points the $Im\{g(s)\} = 0$. (Meissen, 2013). Where
\begin{equation}\label{eq:69}
g(s) = s - a \mbox{ln} s
\end{equation}
in (\ref{eq:68}).
Thus we have  \label{eq:70}
\begin{equation}
g^{'}(s) = 1 - \frac{a}{s}
\end{equation}
So the stationary point occurs when $s = a$.\\
With $s = u+iv$ we have 
\begin{equation}\label{eq:71}
\mbox{Im} \{u + iv -a\mbox{ln}(u + iv)\} = 0.
\end{equation}
Expressing $u + iv$ as $Re^{i\theta }$ we have 
\begin{equation}\label{eq:72}
Im \big\{ (u - a\mbox{ln}R) + i(v - a\theta )\big\} = 0
\end{equation}
then
\begin{equation}\label{eq:73}
v = a\theta 
\end{equation}
and as
\begin{equation}\label{eq:74}
 \theta = arg(s)= \tan^{-1}\bigg(\frac{v}{u}\bigg)
\end{equation}
then 
\begin{equation}\label{eq:75}
u = a \theta \cot\theta 
\end{equation}
(Abate and Valko, 2003).\\

With $s = u + iv$ and $v = a\theta$ it can be shown that  (\ref{eq:75}) can be parametrized to Talbots contour:
\begin{equation}
s(\theta) =a\theta(\cot(\theta) + i),\;\ \;\ -\pi < \theta < +\pi
\end{equation}
(Talbot,1976). 

\subsubsection{Conformal mapping of the Talbot contour.}

While the above parametrization can be used as a basis for inverting the Laplace Transform we proceed with the algorithm's development via a convenient conformal mapping as follows.
\begin{equation}\label{eq:76}
\cot\theta = \frac{i(e^{i\theta} +e^{-i\theta })}{(e^{i\theta} -e^{-i\theta })}
\end{equation}
Then
\begin{equation}\label{eq:77}
\theta \cot\theta  +i\theta 
\end{equation}
\begin{equation}\label{eq:78}
=\frac{2i\theta }{1 - e^{-2i\theta }}
\end{equation}
with $z = 2i\theta $ then (\ref{eq:78}) 
\begin{equation}\label{eq:79}
=\frac{z}{1 - e^{-z}}.
\end{equation}
The function 
\begin{equation}\label{eq:80}
s(z) =\frac{z}{1 - e^{-z}}
\end{equation}
maps the closed interval $M = [-2\pi i,2\pi i]$ on the  imaginary $z-$plane  onto the curve $L$ in the $s$ plane giving the integral,
\begin{equation}\label{eq:81}
f(t) = \frac{1}{2\pi i}\int_{L}F(s)\;\ e^{st}\;\ ds.
\end{equation}
(See Logan,(2001) for the details of this transformation).\\
Next we follow the procedure as adopted by Logan (2001) for numerically integrating (\ref{eq:81}).
With the change of variable (\ref{eq:81}) becomes
\begin{equation}\label{eq:82}
f(t) = \frac{1}{2\pi i} \int_M F(s(z))\;\ e^{s(z)t}\;\ s^{'}(z)\;\ dz
\end{equation}
where 
\begin{equation}\label{eq:83}
s^{'}(z) = \frac{1-(1 + z)e^{-z}}{(1 - e^{-z})^2}.
\end{equation}
For convenience we write,
\begin{equation}\label{eq:84}
f(t) =  \frac{1}{2\pi i}\int_M  I(z,t)\;\ dz
\end{equation}
where 
\begin{equation}\label{eq:85}
I(z,t) = F(s(z))\;\ e^{s(z)t}\;\ s^{'}(z).
\end{equation}
The integral in (\ref {eq:84}) is then rotated by $\frac{\pi}{2}$ so the interval of integration is now real and  becomes $[-2\pi,2\pi]$ and then we use the trapezoid rule with $n$ odd and $w = -iz$ to obtain  
\begin{equation}\label{eq:86}
f(t)\cong \frac{1}{n}\bigg\{(I(2\pi i)+T(-2\pi i) + 2\sum^{n-1}_{j=1}I(iw_j)\bigg\}
\end{equation}
  where 
\begin{equation}
w_j = 2\pi(\frac{2j}{n} -1)
\end{equation}
and we note that $I(2\pi i)= I(-2\pi i) =0$.\\
 Logan (2000).

\subsubsection{The regularization properties of the Talbot algorithm}

Despite the intricacies of deriving the Talbot algorithm we have found it to be a relatively easy algorithm to implement. Also the tests which we have carried out so far show that the algorithm performs to a high degree of accuracy. Moreover the algorithm converges much faster than the Fourier series method without requiring the use of any acceleration schemes. Additionally in the form in which we have used it there is only one parameter to control.\\ 
But perhaps it's greatest strength is the fact that \textit{we have found that it is able to handle noisy data (of magnitude outlined below) with little growth in the corresponding error}. As we will show this is not  the case for either the Fourier series or the Stehfest inversion algorithms presented above. Moreover this ``regularization property" does not exist for many of the numerical inversion schemes as indicated by Egonmwan (2012). For most algorithms this is generally overcome by constructing some regularization scheme which then needs to be  attached on to the inversion algorithm(s) of choice. This of course increases the complexity of the inversion process involving new parameters thus requiring even greater knowledge of the desired solution. This is even more so if the scheme also involves some additional accelerated convergence process.\\
As we pointed out earlier, the perturbation in the numerical schemes are a consequence of the inversion being carried out on the real axis in the complex plane. The inclusion of complex arithmetic in the Talbot scheme enormously diminishes this perturbation.
Of great importance here too is that the \textit{``regularization properties " of the Talbot algorithm means that very good performance can be obtained on many of the test functions without the necessity for multi-precision}.\\
Egonmwan (2012) examines regularized and collocation methods for the numerical inversion of the Laplace transform which involve Tikhonov (1995) based methods.  
 This is then applied to the Stehfest (1970) and Piessens (1972)  methods on various standard test functions for both exact $F(s)$ and noisy $F(s^\delta) $ data.\\ 
 For the Stehfest, Piessens (1972) and the regularized method Egonwan (2012)  added noise of a magnitude $10^{-3}\times rand(1,1)$ to the Laplace transform values. Commenting on his results  Egonwan notes ``the Gaver Stehfest method gave very nice approximate solutions for a wide range of functions. However it completely failed in the presence of noisy data. In the case of exact data the method produced better numerical approximations when compared to the Piessins and the regularized collocation methods. However the Piessins method gave better results than the regularized collocation method in the case of exact data."\\
In other words methods which performed well for exact data did not do well for noisy data and the regularized collocation method failed for exact data. Thus to use such regularized methods requires some a priori knowledge of the magnitude of the noise involved and by implication a better knowledge of the solution than might be otherwise possible.\\
We will demonstrate that the Talbot algorithm performs well when dealing with noisy data. We will also show that both the Fourier series and the Stehfest methods give increasingly bad results for noisy data as we increase the precision in the calculations while the opposite is the case for the Talbot algorithm.\\

\subsection{Tests}

Table 1 lists the functions together with a variety of properties for the purpose of testing the noise handling capability of the three inversion algorithms employed.

\begin{table}[H] 
\centering
\begin{tabular}{|c|c|c|c|}
\hline
{No.}&{$F(s)$}&{$f(t)$}&Function type\\
\hline
1&$\frac{s}{(s^2+1)^2}$&$0.5t\sin (t)$&Oscillating increasing\\
2&$\frac{1}{(s+1)^2}$&$te^{-t}$&Exponentially decreasing\\
3&$\frac{1}{s^5}$&$\frac{1}{24}t^4$&Increasing\\
4&$\frac{1}{\sqrt s}$&$\frac{1}{\pi t}$&With singularities\\
5&$\mbox{erf}\{\frac{2}{\sqrt s}\}$&$\frac{1}{\pi t}\sin (4\sqrt t)$&Oscillating  with singularities\\
6&$\frac{1}{s^2-0.5^2}$&$\sinh (0.5t)$&Hyperbolic\\
7&$\frac{s^3}{s^4+4{(0.5)}^4}$&$\cos(0.5t)\cosh(0.5t)$&Combination of oscillating and hyperbolic\\
&$\frac{\ln s}{s}$&$-(\ln t +\gamma)$&Natural log\\
\hline
\end{tabular}
\caption[Colin's Table]{Test Functions} 
\label{colin}
\end{table}

We use three error measures, the $L^2$ norm defined as
\begin{equation}\label{eq:87}
L^2 = \sqrt{\sum_{i=1}^{40} \bigg|f_{numerical}(t_i)-f_{exact}(t_i)\bigg|^2},\;\ i=1..40
\end{equation}
the $L^{\infty}$ norm as
\begin{equation}\label{eq:88}
L^{\infty} = \max\bigg|f_{numerical}(t_i)-f_{exact}(t_i)\bigg|,\;\ i=1..40
\end{equation}
and the percentage error  as
\begin{equation}
\max\bigg|\frac{f_{numerical}(t_i)-f_{exact}(t_i)}{f_{exact}(t_i)}\times 100\bigg|,\;\ i=1,..40
\end{equation}
With $t$ sampled over $40$ points  for $t=0.1$ to $4$. (These are the same values used by Egonmwan). The $L^2$ norm is chosen as a measure which averages out the error over the sample points while the $L^{\infty}$ norm and the $\%$ error as defined above chooses the maximum error obtained for these measures.
In all cases the Magnitude of noise added is $10^{-3}\times rand(1,1)$)\\
The precision used for implementing the three algorithms is $1.8M$ where $M$ is the number of weights for the Stehfest algorithm and $2N$ where $N$ is the number of terms in the sumation for the Talbot and the Fourier methods. The choice of these levels of precision is based on trial and error from our previous work with these algorithms. \\They are perhaps larger than they need to be but as our interest in this investigation is not on their efficiency but on their ability to handle noisy data we wanted to ensure that the precision played as little part as possible in assessing their performance. Thus in cases where the extended precision decreases the accuracy of the noisy data we used the usual double precision for these inversions.\\
For functions which have sine, cosine and hyperbolic properties we increase the weights for the Stehfest. This is because these functions require more weights and a corresponding increase in precision for the Stehfest method to produce accurate results.
For the Fourier Series method we choose the parameter value of $a=4$. Once again this choice is based on trial and error from our earlier work with this algorithm. We have found that this choice for $a$ gives the best results for inverting the widest class of functions.
\subsection{Results}  
\begin{table}[H] 
\centering 
\begin{tabular}{|c|c|c|c|c|c|c|c|}
\hline
 & &\multicolumn{3}{c|}{$No\;\ Noise$} &\multicolumn{3}{c|}{$Noise$}\\
\hline
 Method&$M$&$L^2$ &$L^\infty$ &$\%error$&$L^2$ &$L^\infty$&$\%error$\\
\hline
Stehfest &30 &9.4(-4) &5.0(-4) &3.8(-2) &4.6(16)&3.6(16)&1.2(18)\\
\hline
Talbot&55 &2.0(-6) &5.4(-7) &2.3(-4) &6.2(-4) &2.7(-4) &3.7(-2)\\
\hline
Fourier&55 &4.2(-2) &1,8(-3) &3.1(-1) &8.9(1) &2.9(0) &1.1(3)\\
\hline
\end{tabular}
\caption[Colin's Table]{ $f(t)=0.5 t \sin (t)$\;\ $= L^{-1} \{\frac{s}{(s^2+1)^2}\}$} 
\label{colin} 
\end{table}
\begin{table}[H] 
\centering 
\begin{tabular}{|c|c|c|c|c|c|c|c|}
\hline
 & &\multicolumn{3}{c|}{$No\;\ Noise$} &\multicolumn{3}{c|}{$Noise$}\\
\hline
 Method&$M$&$L^2$ &$L^\infty$ &$\%error$&$L^2$ &$L^\infty$&$\%error$\\
\hline
Stehfest &16 &1.1(-4) &4.0(-5) &5.4(-1) &3.0(7)&2.4(7)&2.6(10)\\
\hline
Talbot&55 &7.3(-6) &6.4(-6) &2.1(-3) &7.8(-4) &2.3(-4) &3.1(-1)\\
\hline
Fourier&55 &3.6(-3) &1.0(-2) &4.9(-0) &1.1(0) &9.0(-1) &9.7(2)\\
\hline
\end{tabular}
\caption[Colin's Table]{ $f(t)=te^{-t}\;\ =L^{-1}\{\frac{1}{(s+1)^2}\}$} 
\label{colin} 
\end{table}
\begin{table}[H] 
\centering 
\begin{tabular}{|c|c|c|c|c|c|c|c|}
\hline
 & &\multicolumn{3}{c|}{$No\;\ Noise$} &\multicolumn{3}{c|}{$Noise$}\\
\hline
 Method&$M$&$L^2$ &$L^\infty$ &$\%error$&$L^2$ &$L^\infty$&$\%error$\\
\hline
Stehfest &16 &6.7(-6) &3.0(-54) &2.8(-3) &3.8(3)&2.4(3)&1.1(12)\\
\hline
Talbot&55 &3.8(-10) &3.4(-10) &5.1(-4) &2.3(-3) &8.8(-4) &1.5(-1)\\
\hline
Fourier&55 &6.2(-1) &2.9(-1) &2.7(0) &7.6(0) &16.3(1) &2.5(3)\\
\hline
\end{tabular}
\caption[Colin's Table]{ $f(t)=\frac{1}{24}t^4 \;\ =L^{-1}\{\frac{1}{(s)^5}\}$} 
\label{colin} 
\end{table}
\begin{table}[H] 
\centering 
\begin{tabular}{|c|c|c|c|c|c|c|c|}
\hline
 & &\multicolumn{3}{c|}{$No\;\ Noise$} &\multicolumn{3}{c|}{$Noise$}\\
\hline
 Method&$M$&$L^2$ &$L^\infty$ &$\%error$&$L^2$ &$L^\infty$&$\%error$\\
\hline
Stehfest &16 &2.7(-8) &1.3(-8) &7.2(-7) &1.5(7)&1.2(7)&6.5(8)\\
\hline
Talbot&55 &9.2(-2) &9.2(-3) &5.2(-2) &9.2(-2) &9.2(-3) &5.2(-2)\\
\hline
Fourier&55 &6.2(-1) &2.9(-1) &2.7(0) &1.4(1) &6.3(0) &7.1(6)\\
\hline
\end{tabular}
\caption[Colin's Table] {$f(t)=\frac{1}{\sqrt{\pi t}} \;\ =L^{-1}\{\frac{1}{(\sqrt s)}\}$}  
\label{colin} 
\end{table}
\begin{table}[H] 
\centering 
\begin{tabular}{|c|c|c|c|c|c|c|c|}
\hline
 & &\multicolumn{3}{c|}{$No\;\ Noise$} &\multicolumn{3}{c|}{$Noise$}\\
\hline
 Method&$M$&$L^2$ &$L^\infty$ &$\%error$&$L^2$ &$L^\infty$&$\%error$\\
\hline
Stehfest &16 &2.6(-4) &1.6(-4) &6.6(-1) &1.2(7)&9.6(6)&7.2(9)\\
\hline
Talbot&55 &2.2(-2) &2.2(-2) &7.1(-1) &2.2(-1) &2.2(-2) &7.1(-1)\\
\hline
Fourier&55 &1.8(1) &1.1(1) &4.3(3) &3.9(3) &2.2(3) &4.1(6)\\
\hline
\end{tabular}
\caption[Colin's Table] {$f(t)=\frac{1}{\pi t}\sin(4\sqrt t)\;\  =L^{-1} \{erf\big(\frac{2}{\sqrt s}\big)\}$} 
\label{colin} 
\end{table}
\begin{table}[H] 
\centering 
\begin{tabular}{|c|c|c|c|c|c|c|c|}
\hline
 & &\multicolumn{3}{c|}{$No\;\ Noise$} &\multicolumn{3}{c|}{$Noise$}\\
\hline
 Method&$M$&$L^2$ &$L^\infty$ &$\%error$&$L^2$ &$L^\infty$&$\%error$\\
\hline
Stehfest &36 &9.8(-3) &9.2(-3) &2.1(-5) &2.6(7)&2.0(7)&7.0(6)\\
\hline
Talbot&55 &7.2(-6) &7.2(-6) &4.6(-6) &4.5(-4) &3.1(-4) &7.6(-3)\\
\hline
Fourier&55 &1.4(-1) &1.4(-1) &1.9(0) &1.7(1) &5.8(0) &3.4(2)\\
\hline
\end{tabular}
\caption[Colin's Table]{$f(t)=\frac{\sinh(0.5t)}{0.5}\;\  =L^{-1}\bigg\{\frac{1}{s^2-{0.5}^2}\bigg\}$} 
\label{colin} 
\end{table}
\begin{table}[H] 
\centering 
\begin{tabular}{|c|c|c|c|c|c|c|c|}
\hline
 & &\multicolumn{3}{c|}{$No\;\ Noise$} &\multicolumn{3}{c|}{$Noise$}\\
\hline
 Method&$M$&$L^2$ &$L^\infty$ &$\%error$&$L^2$ &$L^\infty$&$\%error$\\
\hline
Stehfest &36/16 &3.7(-4) &3.0(-4) &3.0(-4) &3.1(6)&2.4(6)&1.0(8)\\
\hline
Talbot&55 &5.8(-4) &5.8(-4) &5.8(-1) &7.0(-4) &6.0(-4) &6.0(-2)\\
\hline
Fourier&55 &9.4(-2) &6.0(-2) &3.5(-1) &9.0(1) &2.8(1) &5.2(4)\\
\hline
\end{tabular}
\caption[Colin's Table]{ $f(t)=\cosh(0.5t)\cos (0.5t)=L^{-1}\bigg\{\frac{s^3}{s^4+{0.5}^2}\bigg\}$}
\label{colin} 
\end{table}

\begin{table}[H] 
\centering 
\begin{tabular}{|c|c|c|c|c|c|c|c|}
\hline
 & &\multicolumn{3}{c|}{$No\;\ Noise$} &\multicolumn{3}{c|}{$Noise$}\\
\hline
 Method&$M$&$L^2$ &$L^\infty$ &$\%error$&$L^2$ &$L^\infty$&$\%error$\\
\hline
Stehfest &16 &1.9(-8) &1.2(-7) &2.8(-5) &1.4(7)&1.8(7)&2.4(9)\\
\hline
Talbot&55 &6.9(-3) &6.9(-3) &4.0(-1) &7.1(-3) &7.1(-3) &4.1(-1)\\
\hline
Fourier&55 &8.6(-1) &8.3(-2) &4.0(3) &1.2(2) &3.8(2) &6.3(3)\\
\hline
\end{tabular}

\caption[Colin's Table]{ $f(t)=-(\ln(t) +\gamma)=L^{-1}\bigg\{\frac{\ln s}{s}\bigg\}$}
\label{colin} 
\end{table}
\begin{figure}[H] 
\includegraphics[height=7cm,keepaspectratio=true]{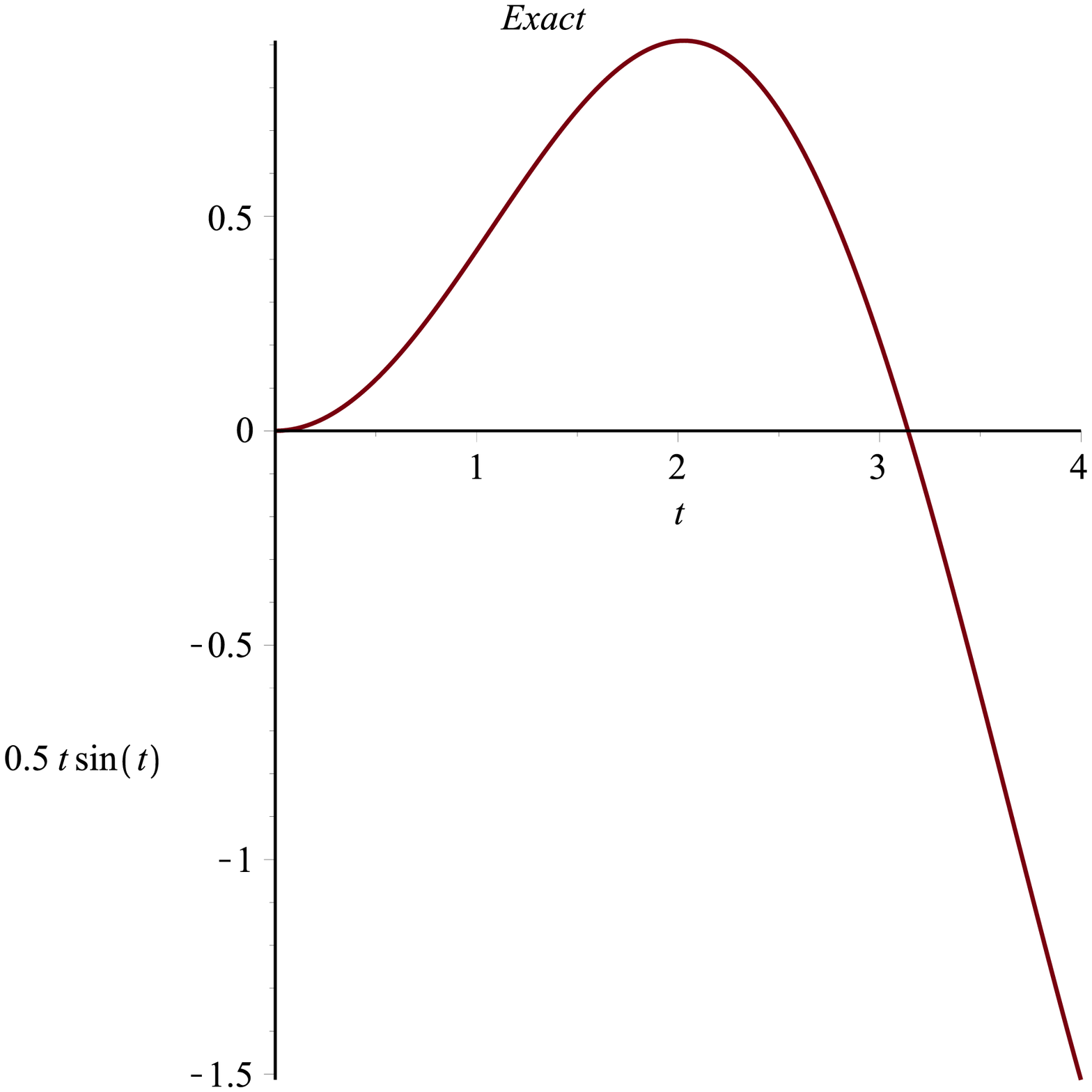}\includegraphics[height=7cm,keepaspectratio=true]{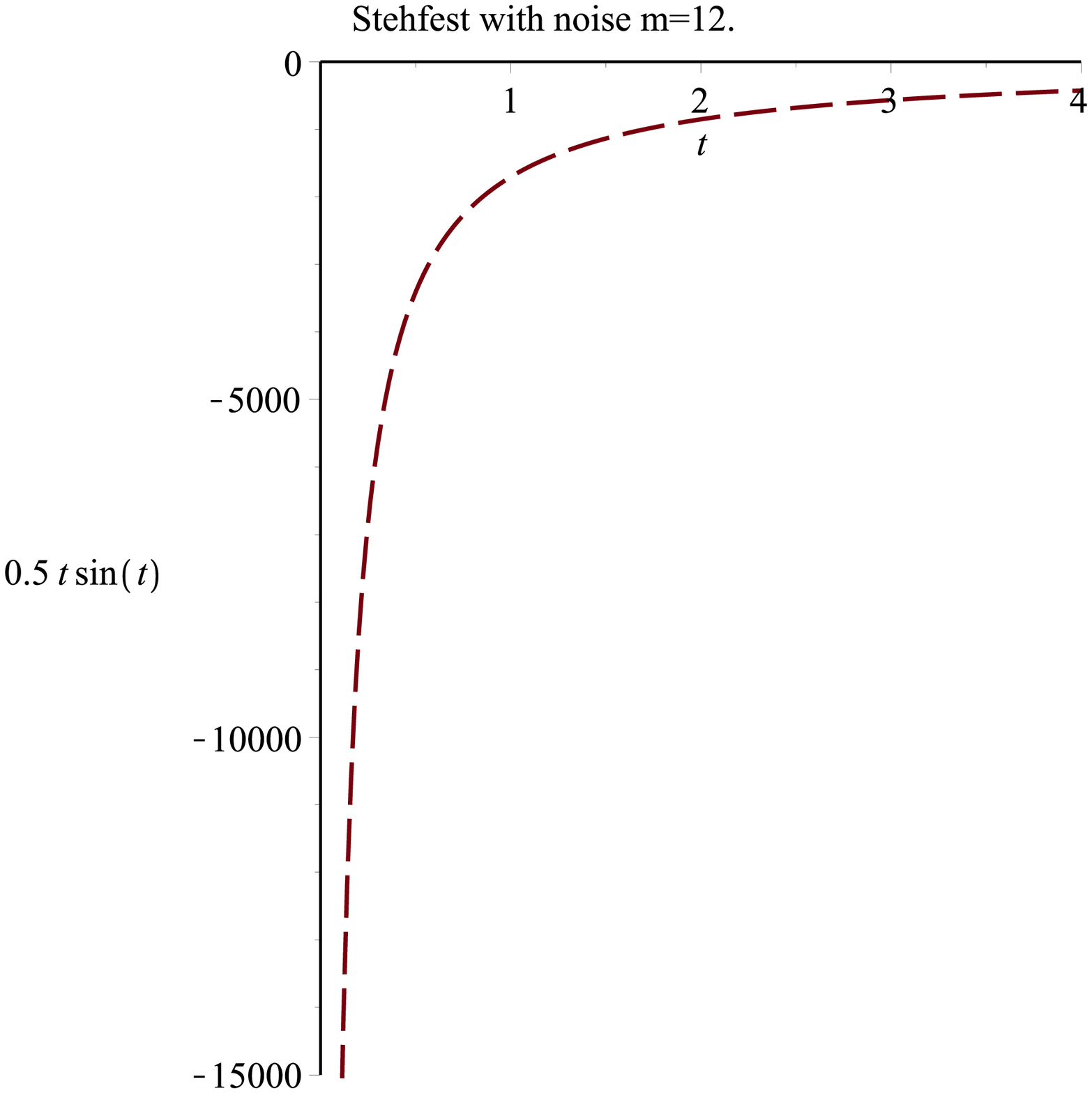}\\
\includegraphics[height=7cm,keepaspectratio=true]{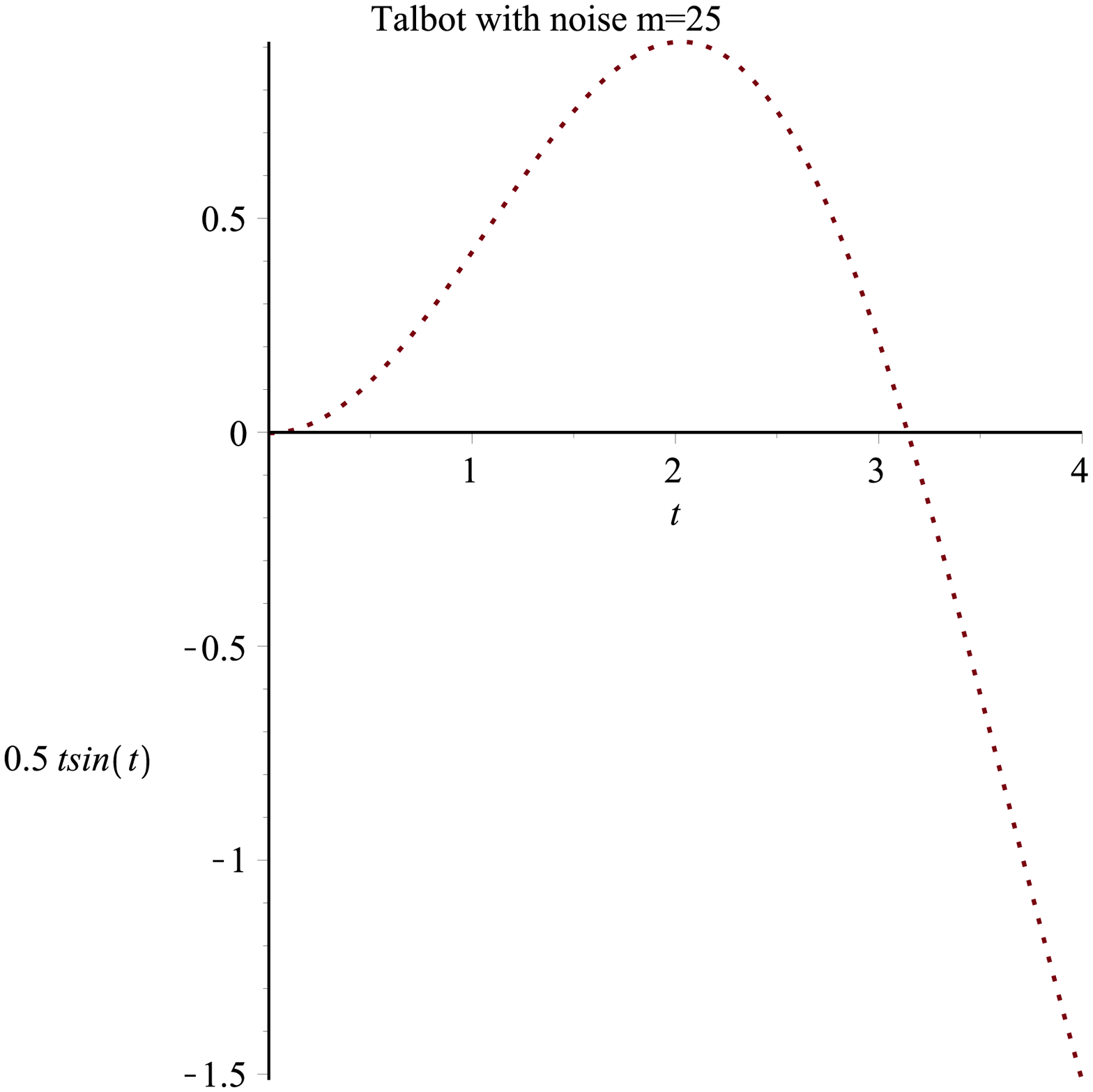} \includegraphics[height=7cm,keepaspectratio=true]{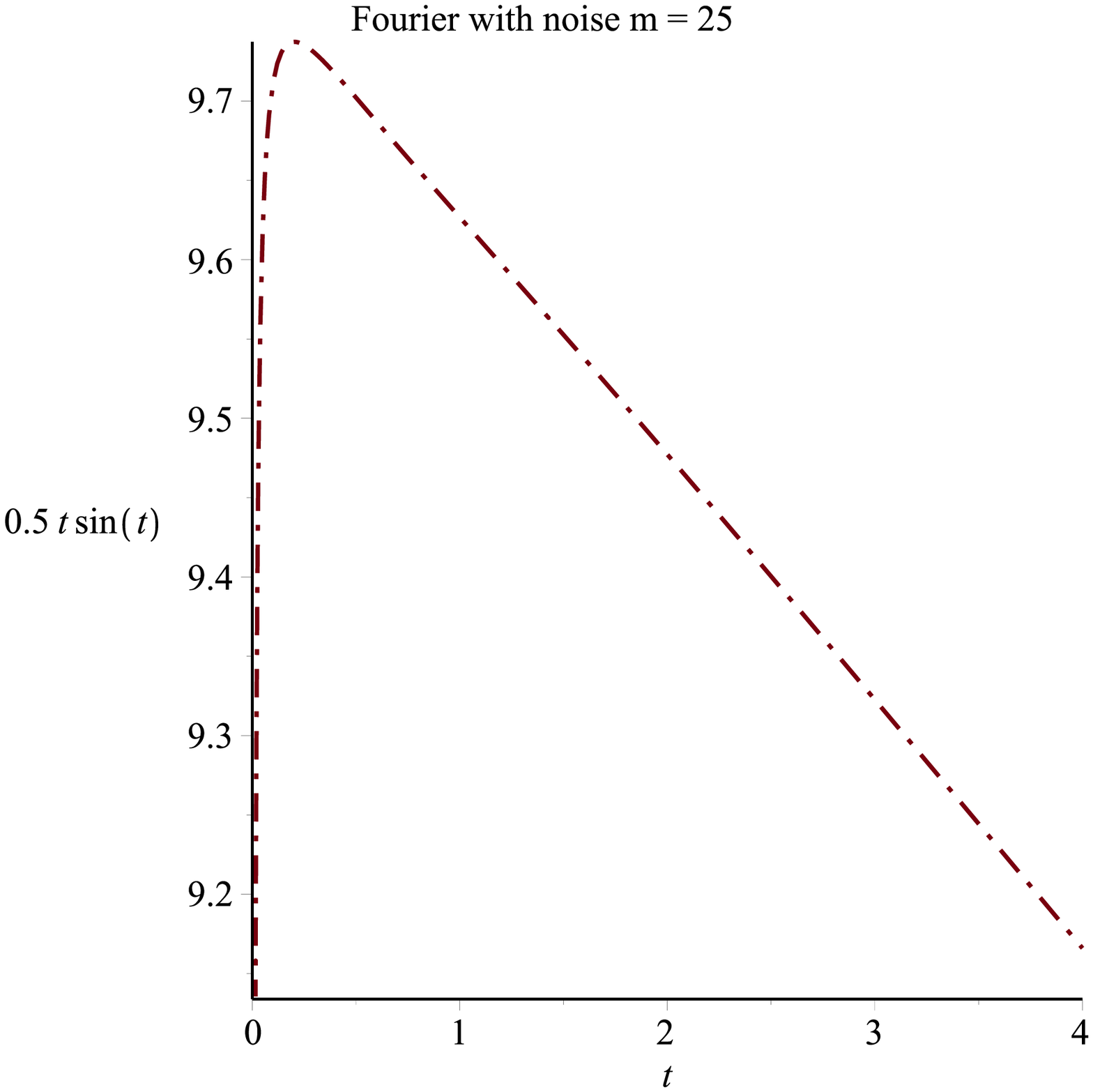}\\ \\ \\

\caption{Numerical Reconstruction of $f(t)=0.5 t. \sin(t)=L^{-1} \{\frac{s}{(s^2+1)^2}\}$}
\end{figure}

\begin{figure}[H] 
\includegraphics[height=7cm,keepaspectratio=true]{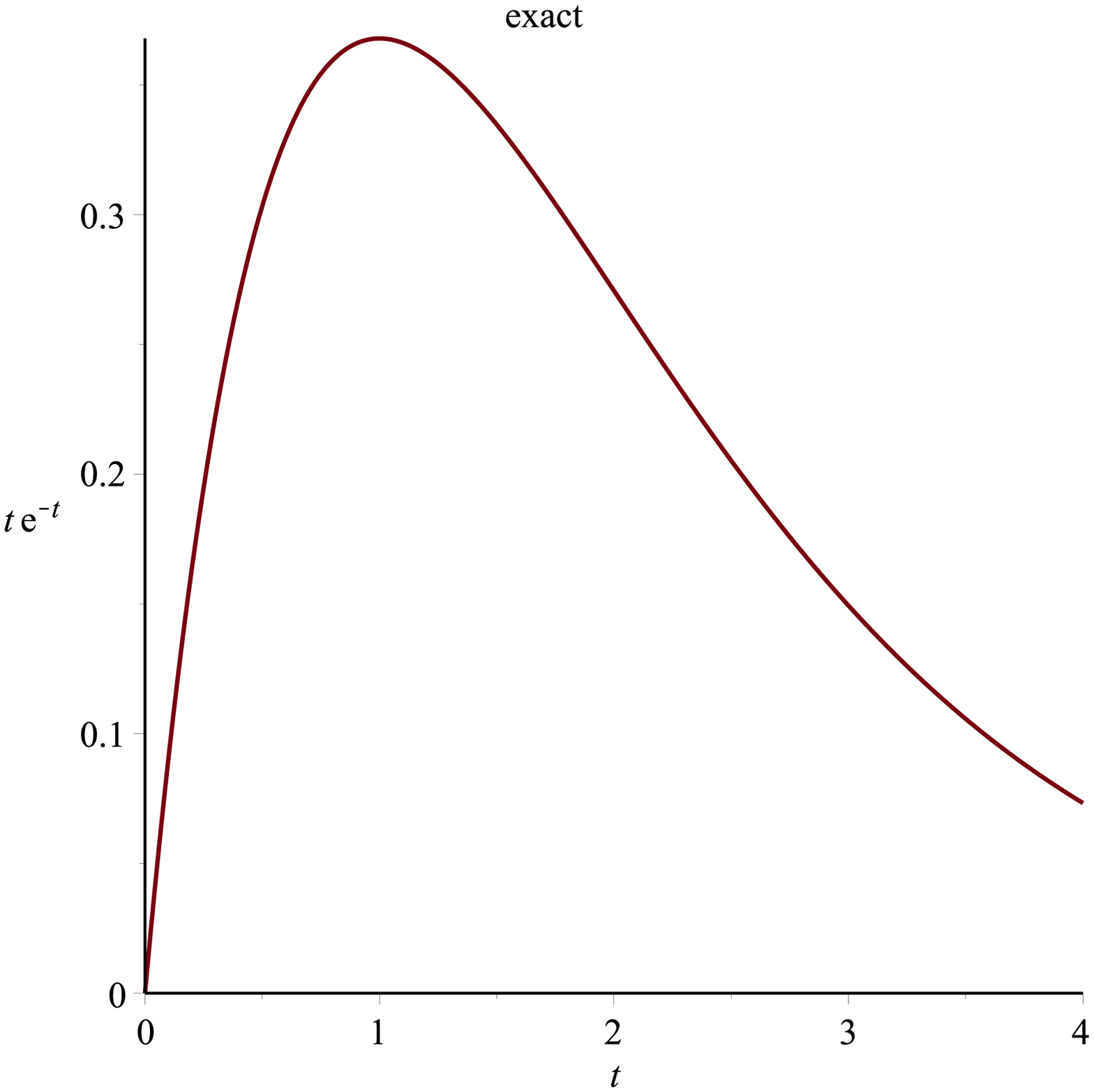} \includegraphics[height=7cm,keepaspectratio=true]{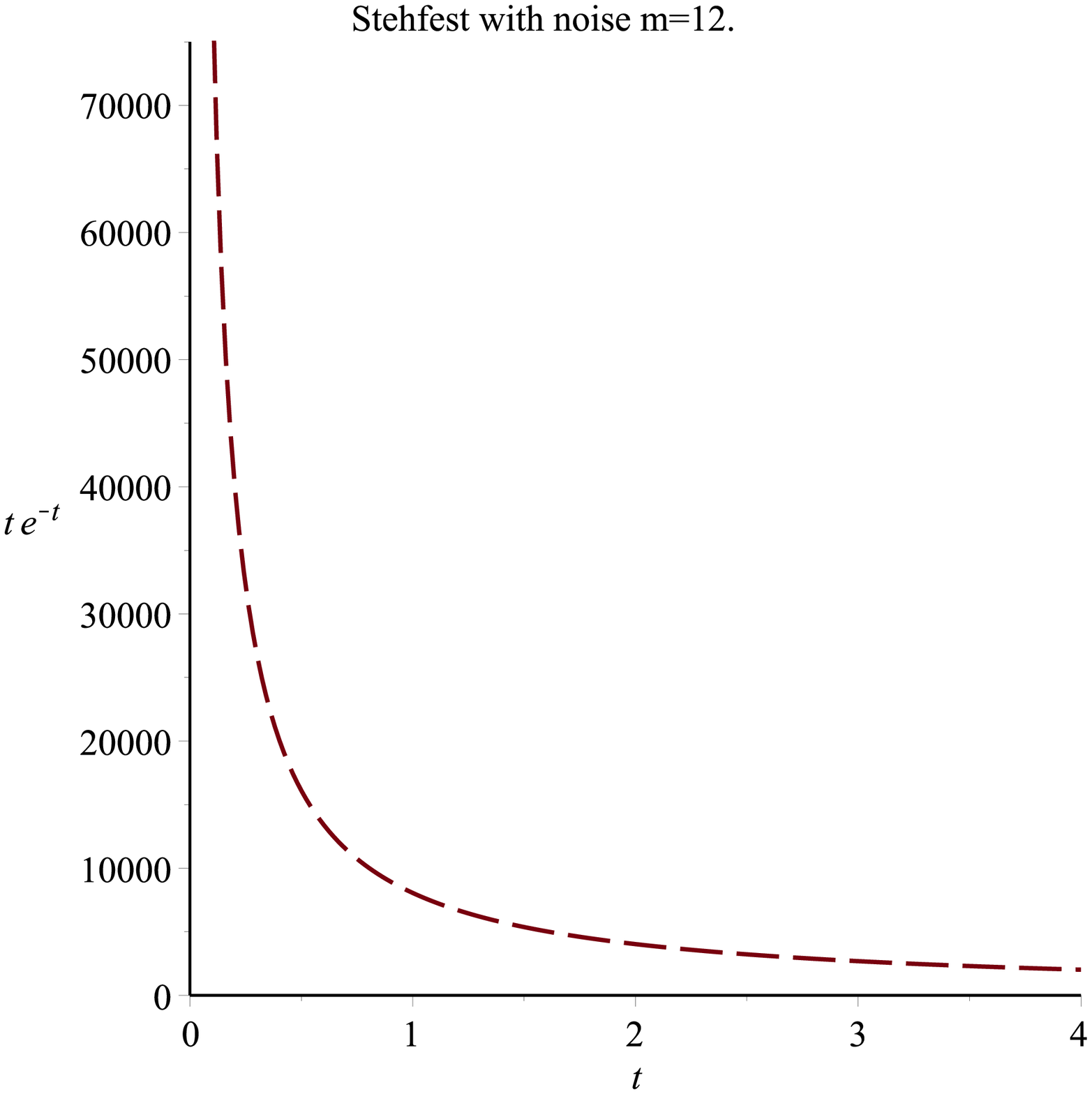}\\
\includegraphics[height=7cm,keepaspectratio=true]{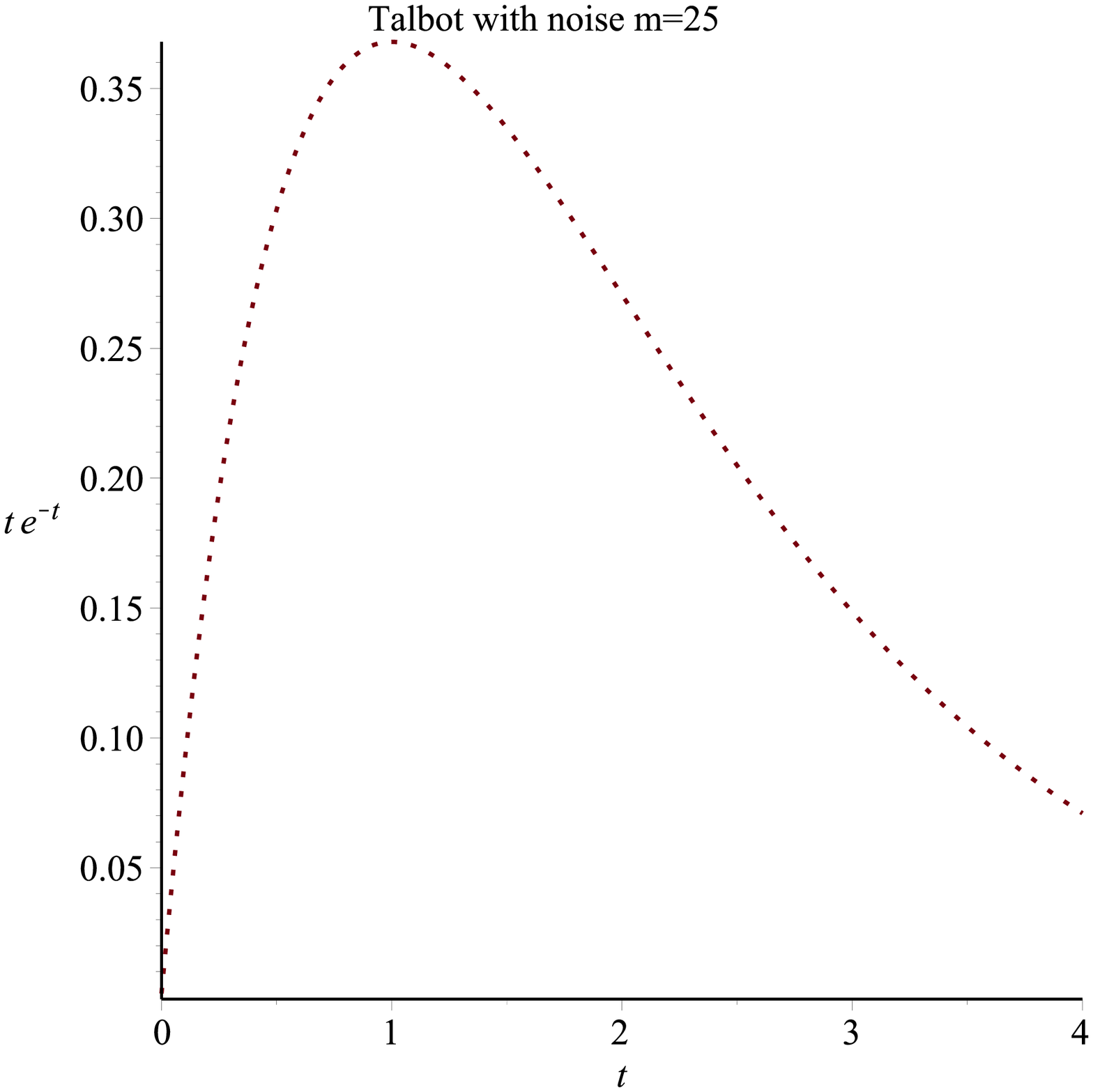}\includegraphics[height=7cm,keepaspectratio=true]{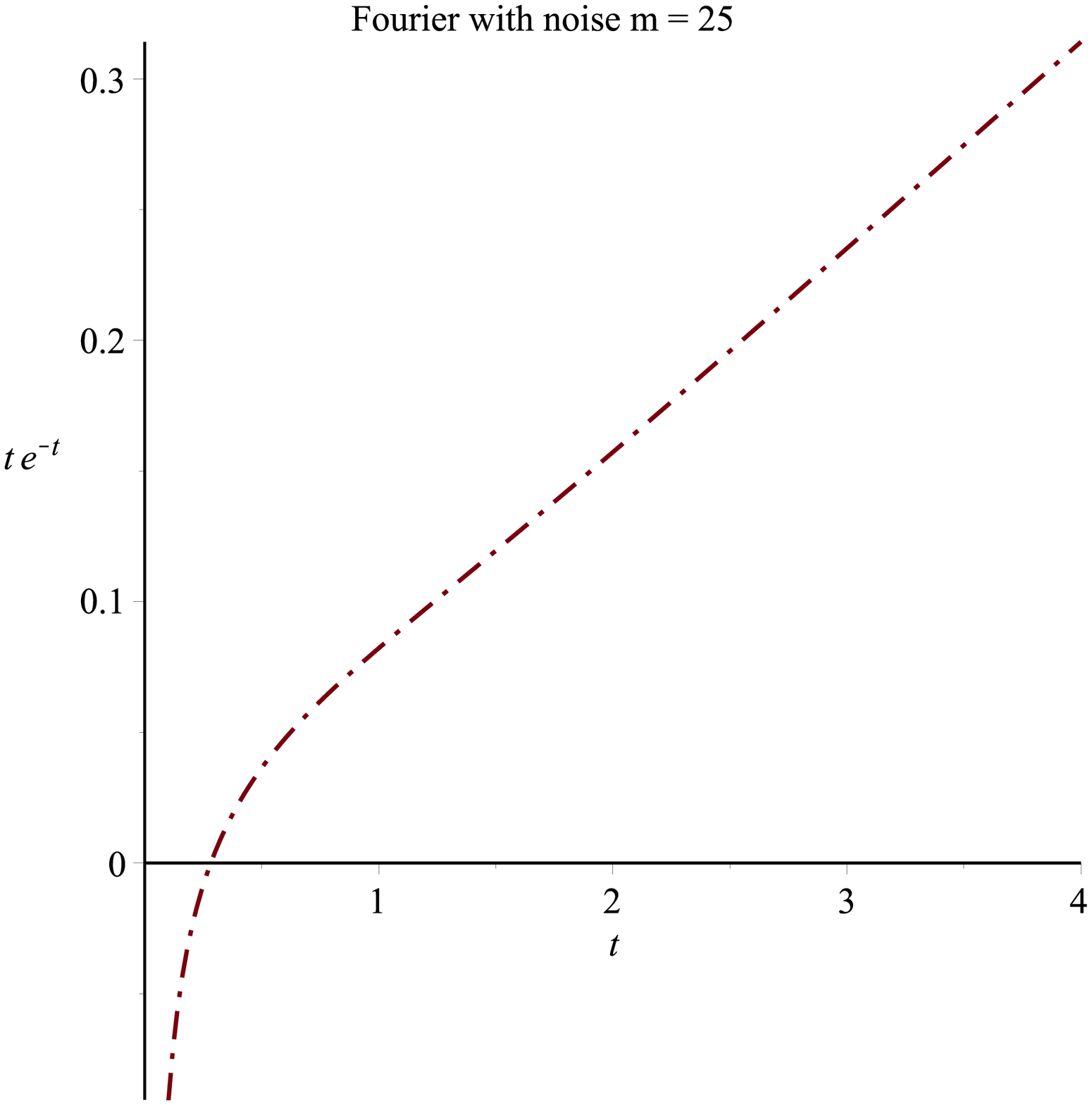}\\ \\ \\

\caption{Numerical Reconstruction of $f(t)= te^{-t} = L^{-1}\{\frac{1}{(s+1)^2}\}$}
\end{figure}
\newpage
\begin{figure}[H] 
\includegraphics[height=7cm,keepaspectratio=true]{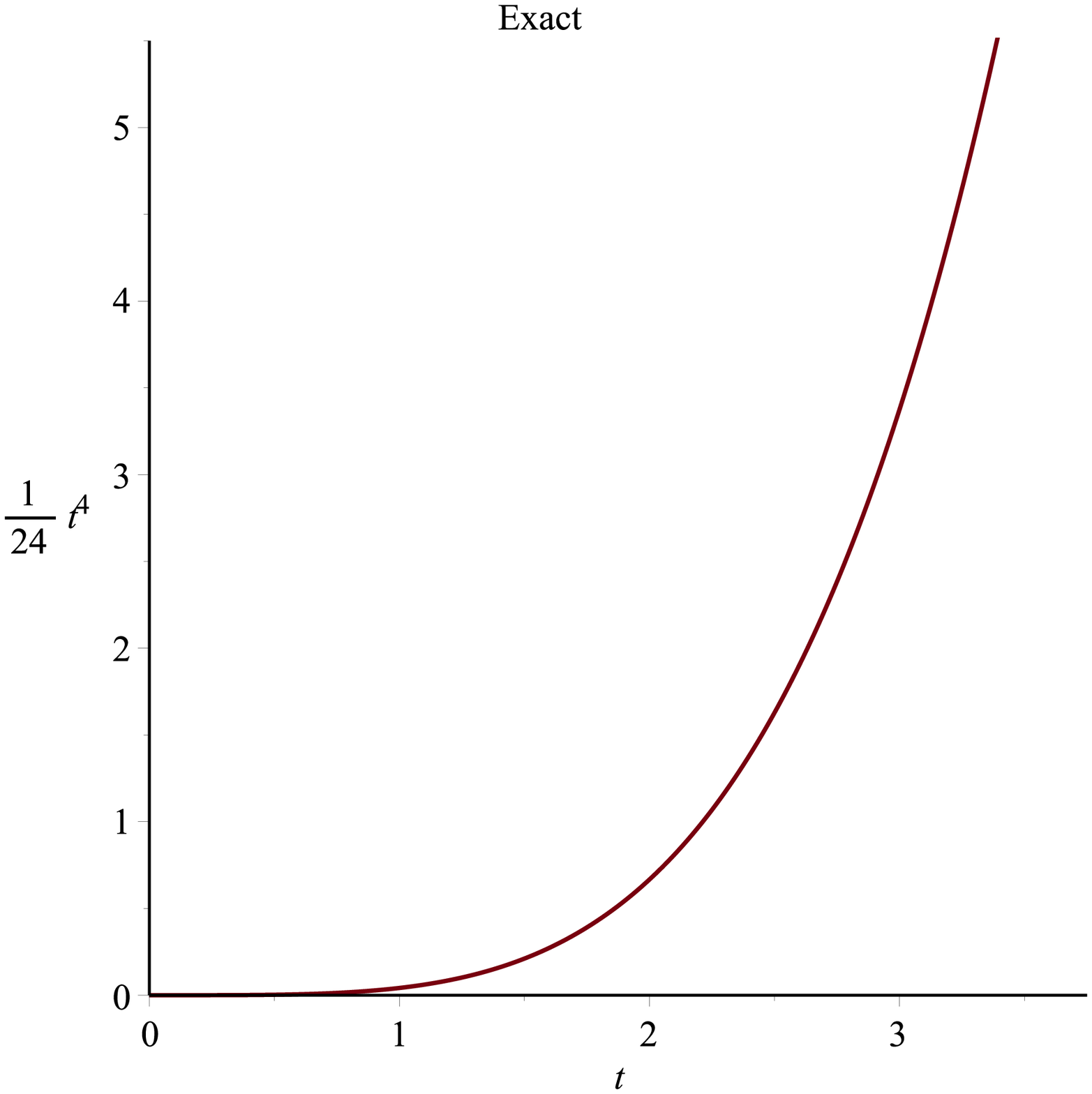} \includegraphics[height=7cm,keepaspectratio=true]{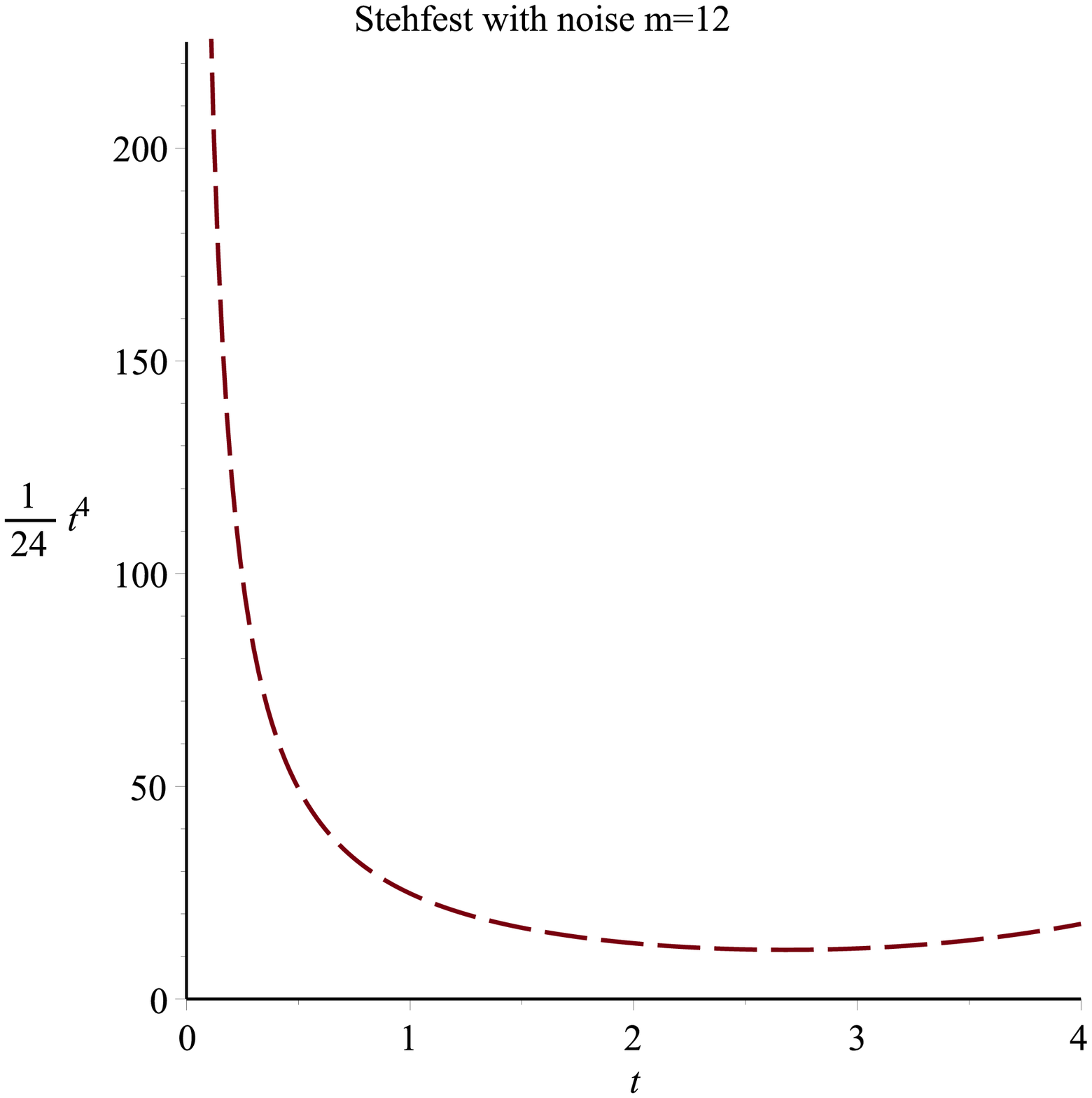}\\
\includegraphics[height=7cm,keepaspectratio=true]{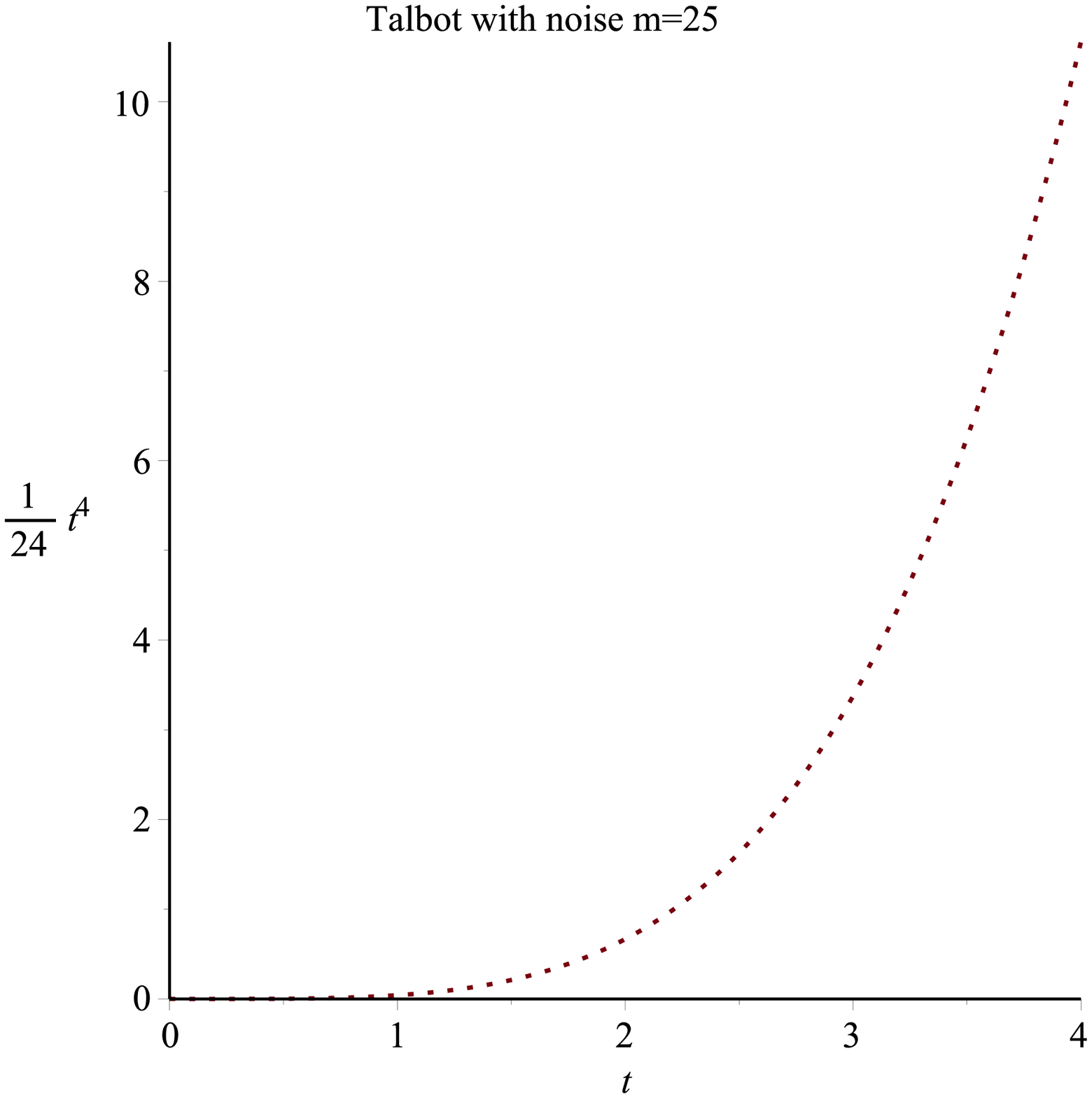}\includegraphics[height=7cm,keepaspectratio=true]{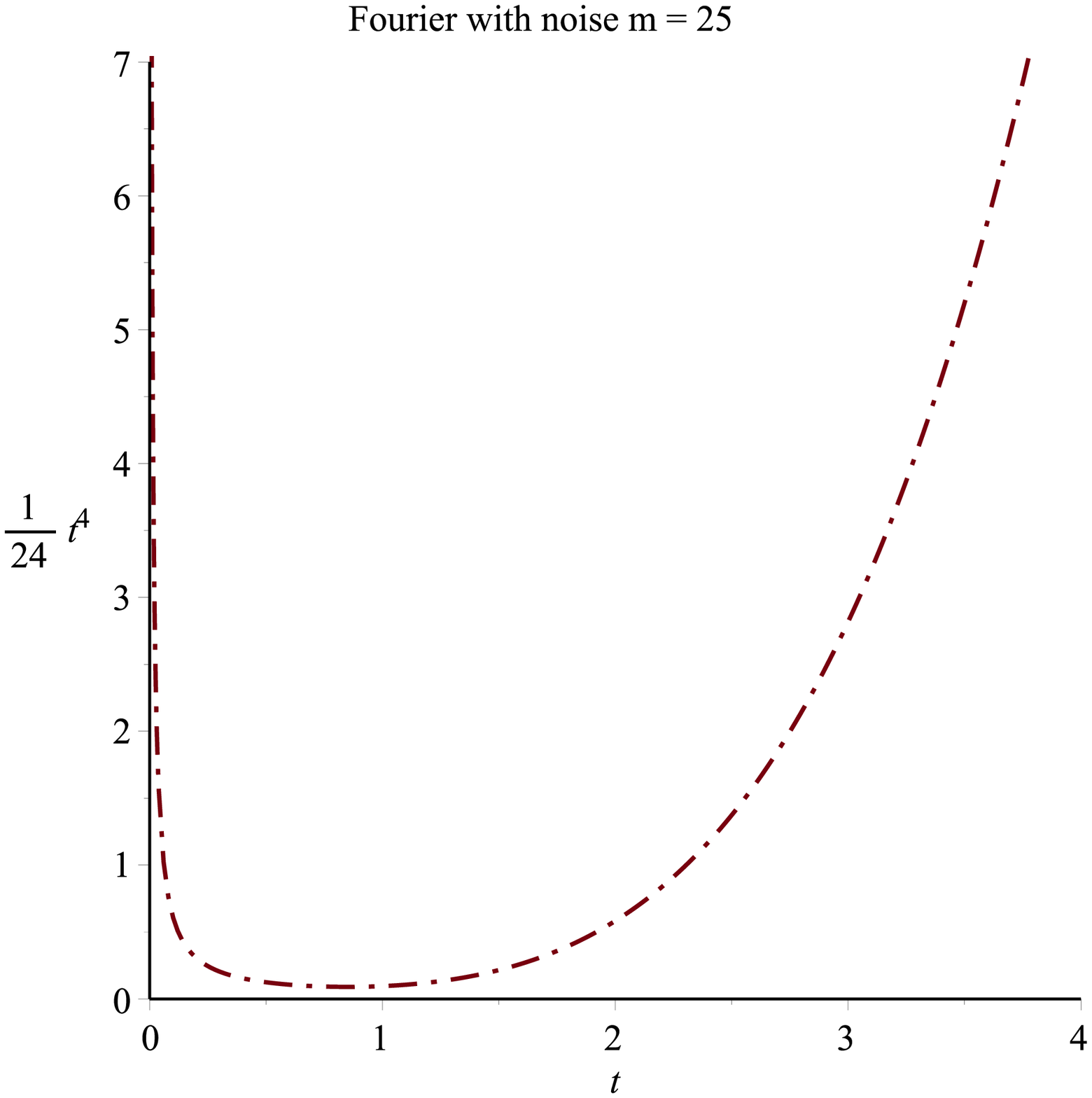}\\ \\ \\

\caption{Numerical Reconstruction of $f(t) = \frac{1}{24}t^4 = L^{-1}\{\frac{1}{s^5}\}$}
\end{figure}
\newpage
\begin{figure}[H] 
\includegraphics[height=7cm,keepaspectratio=true]{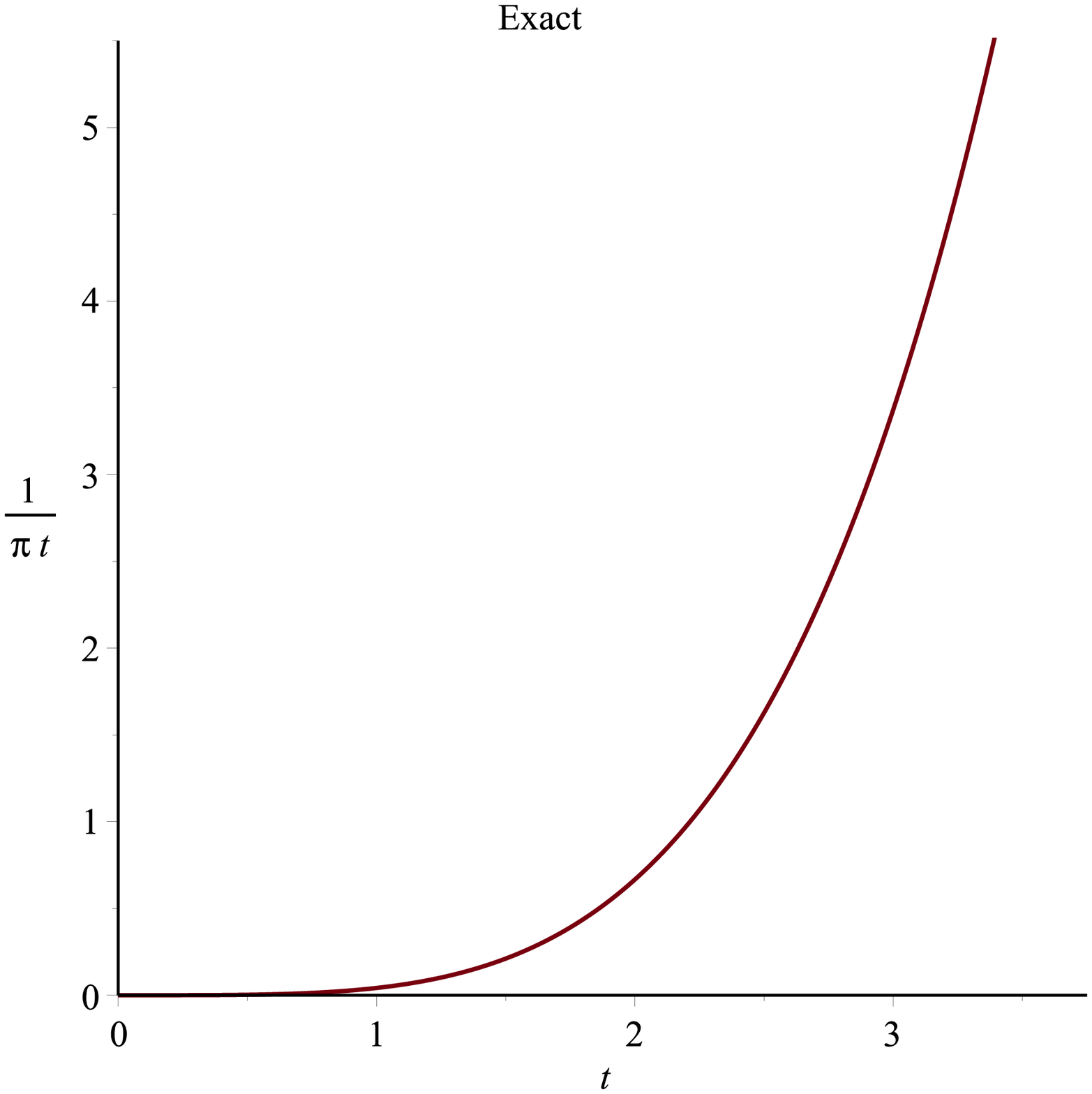} \includegraphics[height=7cm,keepaspectratio=true]{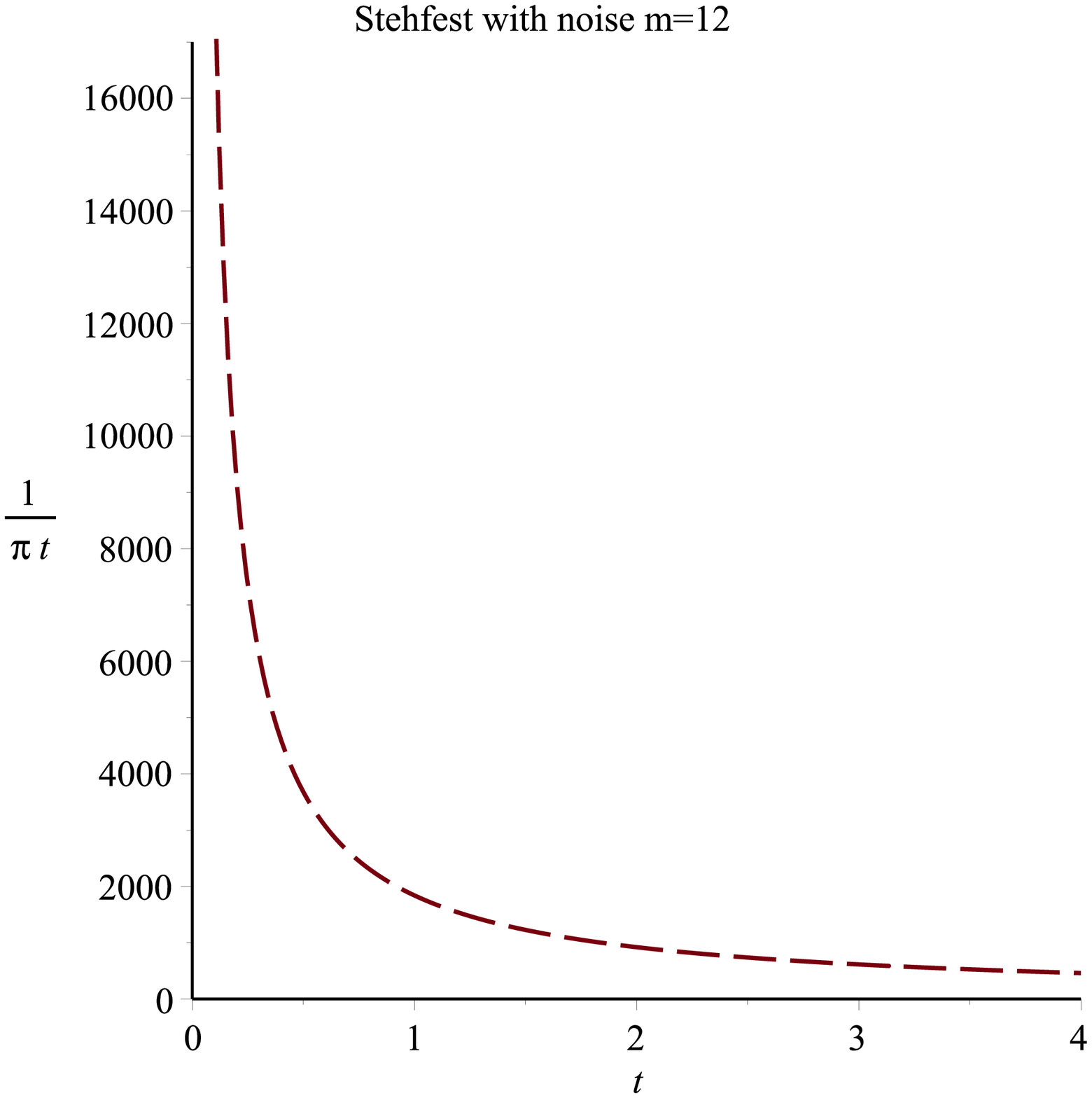}\\
\includegraphics[height=7cm,keepaspectratio=true]{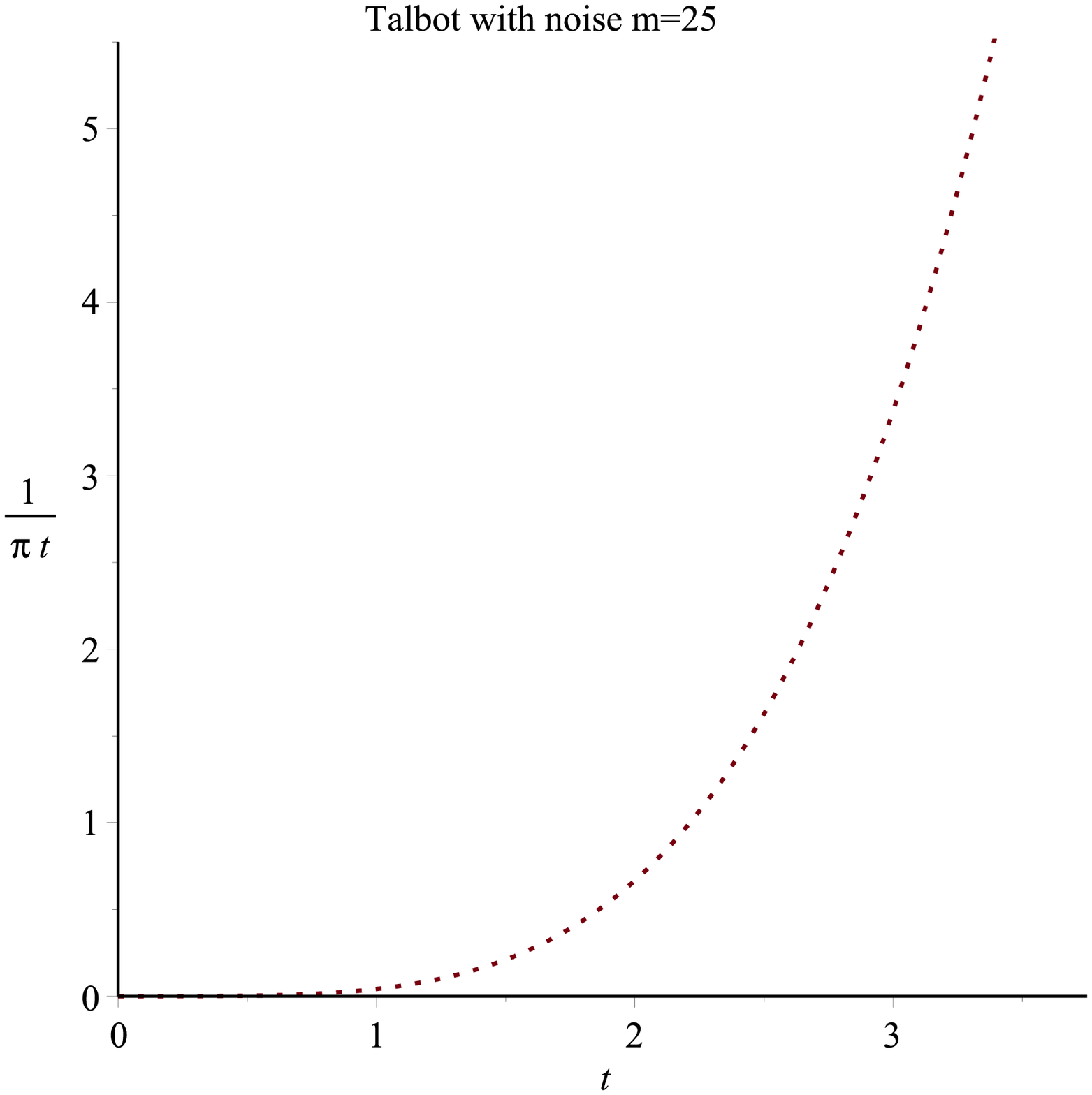}\includegraphics[height=7cm,keepaspectratio=true]{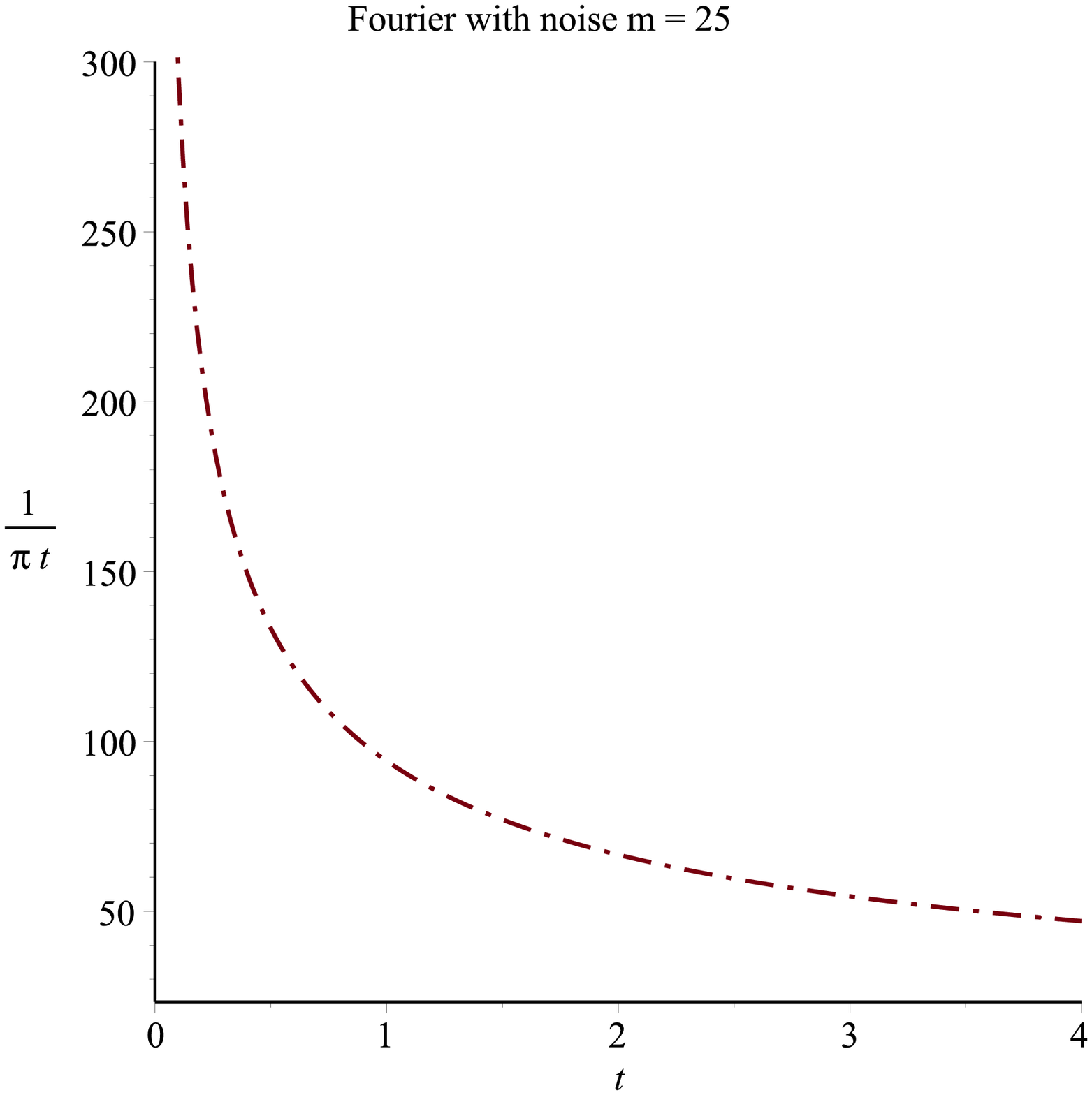}\\ \\ \\

\caption{Numerical Reconstruction of $f(t) = \frac{1}{\pi t} = L^{-1}\{\frac{1}{\sqrt s}\}$}
\end{figure}
\newpage
\begin{figure}[H]  
\includegraphics[height=8cm,keepaspectratio=true]{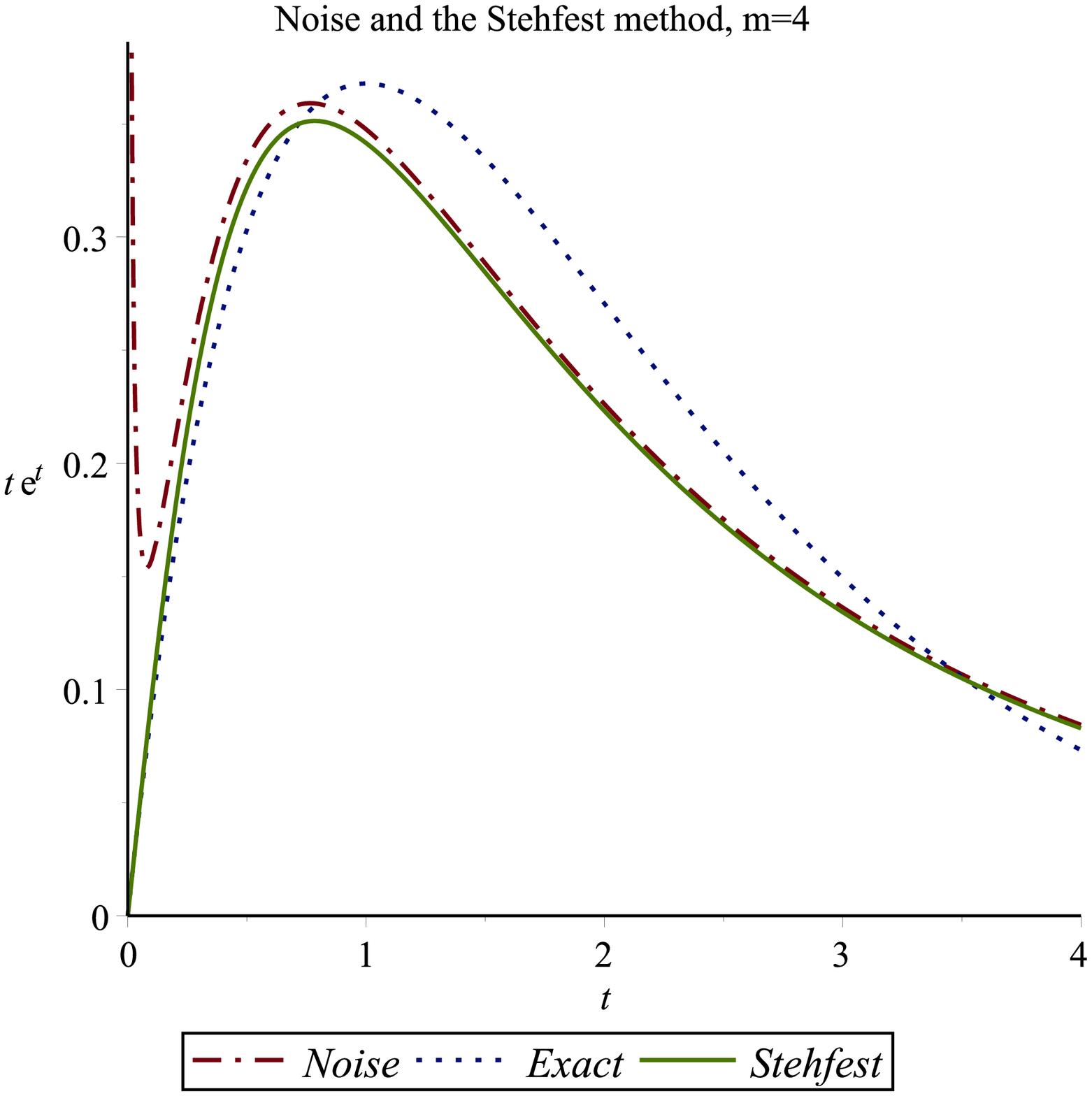} \\ \\ \\
\caption{Numerical Reconstruction of $f(t)=0.5 t. \sin(t)=L^{-1} \{\frac{s}{(s^2+1)^2}\}$}
\end{figure}
\begin{figure}[H]
 
\includegraphics[height=8cm,keepaspectratio=true]{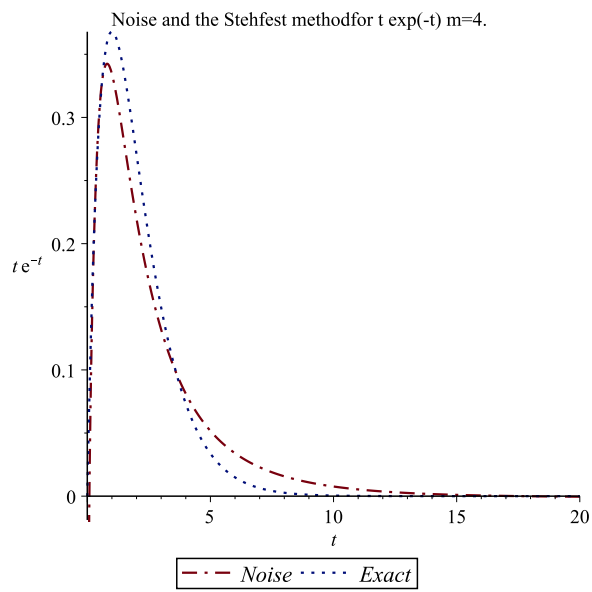} \\ \\ \\

\caption{Numerical Reconstruction of $f(t)=0.5 t. \sin(t)=L^{-1} \{\frac{s}{(s^2+1)^2}\}$}
\end{figure}
\newpage
The above tables (2-9) and graphs (Figs 1 to 4) show very good performance of the Talbot algorithm in handling noisy data. (For brevity we have included only four graphical results for the eight functions using different weights as the performance of these functions with a higher number of weights are well illustrated in the tables). With the exception of the function $f(t) = \frac{1}{\pi t}$ (for which the $L^2$ norm and $L^\infty$ norm  maintain their very small size) the error for the Talbot inversion  diminishes considerably as a function of $M$. However for both the Fourier series and the Stehfest inversion methods both measures of error increase as $M$ increases.\\
The graphs show plots for all three methods and compares them with the exact data. As as can be seen only the Talbot inversion scheme provides an accurate numerical reconstruction of the noisy data. Table 8 includes two sets of weights for the Stehfest inversion algorithm. For the accurate inversion of sinusoidal functions this algorithm requires more weights for increasing values of $t$, here for example we use 36 weights. However when noise is added the accuracy decreases with the number of weight used thus in this case for better performance we have used 16 weights.\\ We observe that the recovery of the function number 5 in table one performs badly for the Fourier series method in both the noisy and noise free environment. Figure 5 and Figure 6 demonstrate that the Stehfest algorithm handles noisy data more accurately by decreasing the number of weights used. This is because the error generated in reconstructing the function from noisy data increases as the number of weights used rises. 
However the accuracy achieved by decreasing the number of weight is not sufficient to justify such an approach for handling noisy data. Moreover as we have stated a larger number of weights and the corresponding increase in precision is necessary for handling trigonometric and hyperbolic functions. We again note that no such considerations are necessary when employing the Talbot algorithm.
\section{Conclusions.}

The results show that the Talbot algorithm handles the noisy data extremely well having very little impact on the final outcome. Both the Stehfest and the Fourier series methods fail to handle the noise. This has implications for implementing the LTFD method when solving nonlinear diffusion or time dependent  parabolic partial differential equations which can generate noisy data through a combination of measurement, truncation and round off error. Using the Talbot algorithm in these circumstances avoids additional complications such as having to devise regularized collocation methods to attain accurate solutions to these problems.

\end{document}